\documentclass[reqno]{amsart}
\usepackage{hyperref}
\newtheorem*{theoA}{Theorem A}
\newtheorem*{theoB}{Theorem B}
\newtheorem*{theoC}{Theorem C}
\newtheorem*{theoD}{Theorem D}
\newtheorem*{theoE}{Theorem E}
\newtheorem*{theoF}{Theorem F}
\newtheorem*{theoG}{Theorem G}
\newtheorem*{theoH}{Theorem H}
\newtheorem*{theoI}{Theorem I}
\newtheorem*{theoJ}{Theorem J}
\newtheorem{theo}{Theorem}[section]
\newtheorem{lem}{Lemma}[section]

\newtheorem{ques}{Question}[section]
\newtheorem{exm}{Example}[section]

\newcommand{\ol}{\overline}
\newcommand{\be}{\begin{equation}}
\newcommand{\ee}{\end{equation}}
\newcommand{\beas}{\begin{eqnarray*}}
	\newcommand{\eeas}{\end{eqnarray*}}
\newcommand{\bea}{\begin{eqnarray}}
\newcommand{\eea}{\end{eqnarray}}

\numberwithin{equation}{section}

\begin{document}
\title[\hfilneg \hfil System of Fermat-type difference equations in $ \mathbb{C}^n $]
{Entire solutions of system of Fermat-type difference and partial differential-difference equations in $ \mathbb{C}^2 $}
\author[G. Haldar \hfil
\hfilneg]
{Goutam Haldar }

\address{Goutam Haldar, Department of Mathematics, Malda College, Malda - 732101, West Bengal, India.}
\email{goutamiit1986@gmail.com, goutamiitm@gmail.com}

\subjclass[2010]{39A45, 30D35, 32H30, 39A14, 35A20.}
\keywords{Several complex variables, meromorphic functions, transcendental entire functions, Fermat-type equations, Nevanlinna theory.}

\maketitle

\begin{abstract}
In this paper we mainly study the existence and the form of entire solutions with finite order for the following system of Fermat-type difference and partial differential-difference equations 
\beas \begin{cases} f_1(z)^2+(\Delta_cf_2(z))^2=1\\
	f_2(z)^2+(\Delta_cf_1(z))^2=1,
\end{cases}\eeas	
\beas\begin{cases} a_1^2f_1(z)^2+(a_2f_2(z+c)+a_3f_2(z))^2=1\\a_1^2f_2(z)^2+(a_2f_1(z+c)+a_3f_1(z))^2=1,\end{cases}\eeas	
\beas\begin{cases} (a_1f_1(z+c)+a_2f_1(z))^2+(a_3f_2(z+c)+a_4f_2(z))^2=1\\(a_1f_2(z+c)+a_2f_2(z))^2+(a_3f_1(z+c)+a_4f_1(z))^2=1,\end{cases}\eeas and 
\beas\begin{cases} (\partial^{I}f_1(z)+\partial^{J}f_1(z))^{n_1}+f_2(z+c)^{m_1}=1\\(\partial^{I}f_2(z)+\partial^{J}f_2(z))^{n_2}+f_1(z+c)^{m_2}=1\end{cases}\eeas
in several complex variables. Some of our results are improvements and extensions of the previous theorems given by Zheng-Xu \cite{Zheng-Xu & Analysis math & 2021}, Xu-Cao \cite{Xu & Cao & 2018}, Xu \textit{et. al.} \cite{Xu-Liu-Li-JMAA-2020} and Li \textit{et. al.} \cite{Li-Zhang-Xu & 2021 & AIMS}. Moreover, we give some examples which are relevant to the content of the paper.	
\end{abstract}

\section{\textbf{Introduction and main results}}
It is well known to all that Nevanlinna theory is an important tool to study value distribution of entire and meromorphic solutions on complex differential equations (see \cite{Hayman & 1964, Laine & 1993, Yi & Yang & 1995}). In 1995, Wiles and Taylor \cite{Tailor & Wiles & 1995, Wiles & Ann. Math. 1995} pointed out that the Fermat-type equation $x^n+y^m=1$, where $m,n\in \mathbb{N}$ does not admit nontrivial solution in rational numbers for $m=n\geq3$, and does exist nontrivial solution in rational numbers for $m=n=2$. Initially, Fermat-type functional equations were investigated by Montel \cite{Montel & Paris & 1927}, Gross \cite{Gross & Bull. Amer. & 1966, Gross & Amer. Math. & 1966}.\vspace{1mm}
\par In 1939, Iyer \cite{Iyer & J. Indian. Math. Soc. & 1939} investigated the solutions of Fermat-type functional equation \beas f(z)^2+g(z)^2=1\eeas and proved that the entire solutions of the above equation are $f(z)=\cos\alpha(z)$ and $g(z)=\sin\alpha(z)$, where $\alpha(z)$ is an entire function, and no other solutions exist.\vspace{1.2mm}
\par In $1970$, Yang \cite{Yang & 1970} considered the following functional equation \bea\label{e1.1} f^n+g^m=1\eea and proved the following interesting result.
\begin{theoA}\cite{Yang & 1970}	There are no non-constant entire solutions of the functional equation $(\ref{e1.1})$, if $m$, $n$ are positive integers satisfying $1/m+1/n<1$.
\end{theoA}
\par After that many researchers started to investigate the existence and the form of entire and meromorphic solutions of some variations of equation \eqref{e1.1} (see \cite{Chen & Gao & CKMS & 2015,Hu & Complex Var & 1995,Laine & 1993,Li Complex Var & 1996,Li & Int. J. & 2004, Liu-Gao & Sci. Math. & 2019,Liu & Cao & EJDE & 2013,Saleeby & Analysis & 1999,Saleeby & Complex Var. & 2004,Yang & Li & 2004}).
\vspace{1.2mm}
\par In recent years, after the development of difference analogues of Nevanlinna theory, specially the development of difference analogous lemma of the logarithmic derivative by Halburd and Korhonen \cite{Halburd & Korhonen & 2006, Halburd & Korhonen & Ann. Acad. & 2006}, and Chiang and Feng \cite{Chiang & Feng & 2008}, independently, many researchers paid their considerable attention to study the existence of entire and meromorphic solutions of complex difference as well as complex differential-difference equations, and obtained a number of important and interesting results in the literature (see \cite{Hal & Kor & lon & 2007,Heit-Kor-Rieppo-Tohge & CMFT & 2001,Latreuch & Mediterr & 2017,Li & Arch. Math. & 2008,Liu-Gao & Sci. Math. & 2019,Liu & Yang & CMFT & 2013,Liu & JMAA & 2009,Liu & Cao & EJDE & 2013,Liu & Cao & Arch. Math. & 2012,Liu & Yang & 2016,Qi-Liu-Yang & CMFT & 2017,Qi-Yang & EJDE & 2013,Rieppo & Acad. Fenn & 2007,Tang & Liao & 2007,Xu-Tu & EJDE & 2016,Xu-Liu-Li & Meditrr & 2019,Zhang & Liao & 2013}).\vspace{1mm}
\par In view of Theorem A, Liu \textit{et. al.} \cite{Liu & Cao & Arch. Math. & 2012} proved that that Fermat-type difference equation $f^n(z)+f^m(z+c)=1$ has no transcendental entire solution when $n>m>1$ or $n=m>2$, and for the case $n=m=2$, the solutions must be of the form  $f(z)=\sin(Az + B)$, where $c(\neq0)$, $B\in\mathbf{C}$ and $A=(4k+1)\pi/2c$, $k$ is an integer. Later, in $ 2013 $, Liu and Yang \cite{Liu & Yang & CMFT & 2013} extended this result by considering the Fermat-type difference equation $f^2(z)+P^2(z)f^2(z+c)=Q(z)$ where $P(z)$ and $Q(z)$ are two non-zero polynomials.\vspace{1.2mm}
\par After that Liu \cite{Liu & JMAA & 2009}, Liu and Dong \cite{Liu-Dong & EDJE & 2015} considered some variations of Fermat-type equations with more general form 
\bea\label{e1.1a} f(z)^2+(f(z+c)-f(z))^2=1,\eea
\bea\label{e1.1b} a_1^2f(z)^2+(a_2f(z+c)+a_3f(z))^2=1,\eea 
\bea\label{e1.1c} (a_1f(z+c)+a_2f(z))^2+(a_3f(z+c)+a_4f(z))^2=1\eea and obtained some results as follows: (i) there is no transcendental entire solutions with finite order of \eqref{e1.1a}. (ii) \eqref{e1.1b} has transcendental entire solutions with finite order if $a_2^2=a_1^2+a_3^2$, and the form of the solution is $f(z)=\cos (az+b)/a_1$. (iii) \eqref{e1.1c} will have finite order transcendental entire solutions if $a_1^2+a_3^2=a_2^2+a_4^2$, and the solution will be of the form $f(z)=(a_3\cos(aiz+bi)+a_1\sin(aiz+bi))/(a_2a_3-a_1a_4)$, where $a_j$'s are non-zero constants in $\mathbb{C}$ with $a_2a_3-a_1a_4\neq0$.\vspace{1mm} \par Hereafter, we denote $ z+w=(z_1+w_1, z_2+w_2) $ for any $ z=(z_1, z_2) $, $ w=(w_1, w_2) $ and $ c =(c_1, c_2) $, where $z,w,c\in\mathbb{C}^2$ except otherwise stated.\vspace{1mm}
\par In 2018,  Xu and Cao \cite{Xu & Cao & 2018} extended Theorem 1.1 of Liu \textit{et. al.} \cite{Liu & Cao & Arch. Math. & 2012} to several complex variables as follows. 
\begin{theoB}\cite{Xu & Cao & 2018}
Let $c=(c_1,c_2,\ldots,c_n)\in \mathbb{C}^n\setminus\{0\}$. Then any non-constant entire solution $f:\mathbb{C}^n\rightarrow\mathbb{P}^{1}(\mathbb{C})$with finite order of the Fermat-type difference equation $f(z)^2+f(z+c)^2=1$ has the form
of $f(z)=\cos(L(z)+B)$, where $L$ is a linear function of the form $L(z)=a_1z_1+\cdots+ a_nz_n$ on $\mathbb{C}^n$ such that $L(c)=-\pi/2-2k\pi$ $(k\in\mathbb{Z})$, and $B$ is a constant on $\mathbb{C}$.
\end{theoB}
Considering equations \eqref{e1.1a}--\eqref{e1.1c}, Zheng and Xu \cite{Zheng-Xu & Analysis math & 2021}, in 2021, extended the results due to Liu \cite{Liu & JMAA & 2009}, Liu and Dong \cite{Liu-Dong & EDJE & 2015} to the case of several complex variables and obtained the results as follows. 
\begin{theoC}\cite{Zheng-Xu & Analysis math & 2021}
Let $c=(c_1,c_2)\in\mathbb{C}^2\setminus\{0\}$. Then there are no transcendental entire solutions $f:\mathbb{C}^2\rightarrow\mathbb{P}^{1}(\mathbb{C})$ with finite order of equation \eqref{e1.1a}.
\end{theoC}
\begin{theoD}\cite{Zheng-Xu & Analysis math & 2021}
Let $c=(c_1,c_2)\in\mathbb{C}^2\setminus\{0\}$ and $a_1,a_2,a_3$ be nonzero constants in $\mathbb{C}$. If the equation \eqref{e1.1b} has a transcendental entire solution $f:\mathbb{C}^2\rightarrow\mathbb{P}^{1}(\mathbb{C})$ with finite order, then $a_1^2+a_3^2= a_2^2$ and $f(z)$ is of the form
\beas f(z)=\frac{1}{a_1}\sin(L(z)+\Phi(t)+A),\eeas where $L(z)= \alpha_1z_1+ \alpha_2z_2$, $\alpha_1, \alpha_2, A\in\mathbb{C}$, $\Phi(t)$ is a polynomial in $t:=c_2z_1-c_1z_2$ in $\mathbb{C}$, and $L(z)$ satisfies \beas L(c)=\alpha_1c_1+\alpha_2c_2=\theta+k\pi \pm \frac{\pi}{2},\;\;\; \tan \theta=\frac{a_3}{a_1}.\eeas
\end{theoD}
\begin{theoE}\cite{Zheng-Xu & Analysis math & 2021}
Let $c=(c_1,c_2)\in\mathbb{C}^2\setminus\{0\}$, $a_1,a_2,a_3,a_4$ be nonzero constants in $\mathbb{C}$, and let $D :=a_1a_4-a_2a_3\neq$0. If equation \eqref{e1.1c} has a transcendental entire solution $f:\mathbb{C}^2\rightarrow\mathbb{P}^{1}(\mathbb{C})$ with finite order, then $a_1^2+a_3^2=a_2^2+a_4^2$ and $f(z)$ is of the form
\beas f(z)=\frac{1}{2D}\left[-(a_3+ia_1)e^{L(z)+\Phi(t)+A}-(a_3-ia_1)e^{-(L(z)+\Phi(t)+A)}\right],\eeas
where $L(z)= \alpha_1z_1+ \alpha_2z_2$, $\alpha_1, \alpha_2, A\in\mathbb{C}$, $\Phi(t)$ is a polynomial in $t:=c_2z_1-c_1z_2$ in $\mathbb{C}$, and $L(z)$ satisfies \beas e^{L(c)}=e^{\alpha_1c_1+\alpha_2c_2}=-\frac{a_3-ia_1}{a_4-ia_2}=-\frac{a_4+ia_2}{a_3+ia_1}.\eeas
\end{theoE}
\par Now, we consider some system of Fermat-type functional equations as follows.\vspace{1mm}

\bea\label{e1.2} \begin{cases} f_1(z_1,z_2)^2+(\Delta_cf_2(z_1,z_2))^2=1\\
	f_2(z_1,z_2)^2+(\Delta_cf_1(z_1,z_2))^2=1,
\end{cases}\eea	where $c=(c_1,c_2)$ be a constant in $\mathbb{C}^2$.
\bea\label{e1.3a}\begin{cases} a_1^2f_1(z)^2+(a_2f_2(z+c)+a_3f_2(z))^2=1\\a_1^2f_2(z)^2+(a_2f_1(z+c)+a_3f_1(z))^2=1,\end{cases}\eea
\bea\label{e1.3b}\begin{cases} (a_1f_1(z+c)+a_2f_1(z))^2+(a_3f_2(z+c)+a_4f_2(z))^2=1\\(a_1f_2(z+c)+a_2f_2(z))^2+(a_3f_1(z+c)+a_4f_1(z))^2=1,\end{cases}\eea where $f_j:\mathbb{C}^2\rightarrow\mathbb{P}^{1}(\mathbb{C})$, $j=1,2$, $c=(c_1,c_2)$ be a constant in $\mathbb{C}^2\setminus\{0\}$, $a_1,a_2,a_3,a_4$ are nonzero constants in $\mathbb{C}$ and $\Delta_cf(z)=f(z_1+c_1,z_2+c_2)-f(z_1,z_2)$ as defined in \cite{Korhonen & CMFT & 2012}.\vspace{1.2mm}
\par As far as our knowledge is concerned, although there are some important and remarkable results about the existence and forms of transcendental entire solutions of Fermat-type difference and partial differential-difference equations in several complex variables (see \cite{Hu & Yang & 2014,Xu & Cao & 2018,Xu & Cao & 2020,Xu & wang & 2020,Xu-Meng-Wan & ADE & 2021,Zheng-Xu & Analysis math & 2021}), there are only few results about the solutions of the system of Fermat-type equations in  the literature (see \cite{Gao & Acta Math. Sinica & 2016,Liu-Xu & JM & 2021,Xu-Liu-Li-JMAA-2020}). Some of these results are listed as follows.\vspace{1.2mm}
\begin{theoF}\cite{Xu-Liu-Li-JMAA-2020}
Let $c=(c_1,c_2)$ be a constant in $\mathbb{C}^2$. Then any pair of transcendental entire solutions with finite order for the system of Fermat-type difference equations \beas \begin{cases} f_1(z_1,z_2)^2+(f_2(z_1+c_1,z_2+c_2))^2=1\\
	f_2(z_1,z_2)^2+(f_1(z_1+c_1,z_2+c_2))^2=1,
\end{cases}\eeas have the following forms \beas (f_1(z),f_2(z))=\left(\frac{e^{L(z)+B_1}+e^{-(L(z)+B_1)}}{2},\frac{A_{21}e^{L(z)+B_1}+A_{22}e^{-(L(z)+B_1)}}{2}\right),\eeas where $L(z)=\alpha_1z_1+\alpha_2z_2$, $B_1$ is a constant in $\mathbb{C}$, and $c, A_{21}, A_{22}$ satisfy one of the following cases \begin{enumerate}
\item [(i)] $L(c)=2k\pi i$, $A_{21}=-i$ and $A_{22}=i$, or $L(c)=(2k+1)\pi i$, $A_{21}=i$ and $A_{22}=-i$, here and below $k$ is an integer;\\
\item[(ii)] $L(c)=(2k+1/2)\pi i$, $A_{21}=-1$ and $A_{22}=-1$, or $L(c)=(2k-1/2)\pi i$, $A_{21}=1$ and $A_{22}=1$.
\end{enumerate}
\end{theoF}

\par Motivated by Theorems B--E, one may ask the following question.
\begin{ques}
What can be said about the existence and the forms of transcendental entire solutions with finite order for the system of the Fermat-type functional equations \eqref{e1.2}--\eqref{e1.3b}?
\end{ques}
\par The main purpose of this paper is to investigate the existence and form of transcendental entire solutions with finite order of system of nonlinear Fermat-type functional equations \eqref{e1.2}--\eqref{e1.3b} with the help of Nevanlinna theory and difference logarithmic lemma in several complex variables (see \cite{Cao & Korhonen & 2016,Korhonen & CMFT & 2012}). We extends Theorems B--E from the complex Fermat-type difference equations to the Fermat-type system of difference equations. Here we list our main results as follows.
\begin{theo}\label{t1}	There is no pair of transcendental entire solutions with finite order of the system of Fermat-type difference equation \eqref{e1.2}.
\end{theo}
\begin{theo}\label{t2}
Let $a_1,a_2,a_3$ be three non-zero complex constants in one variable and $c=(c_1,c_2)\in \mathbb{C}^2\setminus\{0\}$. If $(f_1,f_2)$ is a pair of transcendental entire solution with finite order of simultaneous Fermat-type difference equation \eqref{e1.3a}, then $(f_1,f_2)$ takes one of the following form
\begin{enumerate}
\item[I.] $ (f_1(z),f_2(z))=\displaystyle\left(\frac{1}{a_1}\cos(L(z)+\Phi(t)+A),\displaystyle\frac{1}{a_1}\cos(L(z)+\Phi(t)+A+k)\right)$, where $a_2^2=a_1^2+a_3^2$, $e^{2ik}=1$, $e^{2iL(c)}=-\displaystyle\frac{a_1-ia_3e^{-ik}}{a_1+ia_3e^{ik}}$, $L(z)=\alpha_1 z_1+\alpha_2 z_2$ with $A,\;\alpha_1,\; \alpha_2, k\in\mathbb{C}$ and $\Phi(t)$ is a polynomial in $t:=c_2z_1-c_1z_2$. 
\item[II.]  $(f_1(z),f_2(z))=\displaystyle\left(\frac{1}{a_1}\cos(-(L(z)+\Phi(t)+A)+k),\displaystyle\frac{1}{a_1}\cos(L(z)+\Phi(t)+A)\right)$, where $L(z)$, $\Phi(t)$, $A$ are defined as in $I$, satisfying one of the following conditions:\vspace{1mm}\begin{enumerate}
\item [(a)] $e^{iL(c)}=1$, $e^{ik}=\pm i$ and $a_1=\pm(a_2+a_3)$;\vspace{1mm}
\item [(b)] $e^{iL(c)}=-1$, $e^{ik}=\pm i$ and $a_1=\pm (a_2-a_3)$;\vspace{1mm}
\end{enumerate} 
\item[III.]  $(f_1(z),f_2(z))=\left(\displaystyle\frac{\cos(-(L(z)+\Phi(t)+A)+L(c)+k)}{a_1},\displaystyle\frac{\cos(L(z)+\Phi(t)+A)}{a_1}\right)$, where $L(z)$, $\Phi(t)$, $A$ are defined as in $I$, satisfying one of the following conditions:\vspace{1mm}\begin{enumerate}	\item [(a)] $e^{iL(c)}=1$, $e^{ik}=\pm i$ and $a_1=\pm(a_2+a_3)$;\vspace{1mm}	
\item [(b)] $e^{iL(c)}=-1$, $e^{ik}=\pm i$ and $a_1=\pm (a_2-a_3)$;\vspace{1mm}
\end{enumerate} 
\end{enumerate}
\end{theo}
\par The following examples show the existence of transcendental entire solutions with finite order of the system \eqref{e1.3a}.
\begin{exm}
Let $a_1=3$, $a_2=5$, $a_3=4$ and $L(z)=7z_1-5z_2$. Choose $k\in\mathbb{C}$ such that $e^{ik}=1$. Also, let $c=(c_1,c_2)\in\mathbb{C}^2$ such that $e^{iL(c)}=\cos\frac{2\pi+\alpha}{2}+i\sin\frac{2\pi+\alpha}{2}$, where $\tan\alpha=\frac{24}{7}$. Then, it can be easily verified that \beas &&(f_1(z),f_2(z))\\&&=\left(\frac{\cos(7z_1-5z_2+i(c_2z_1-c_1z_2)^3+3)}{3},\frac{\cos(7z_1-5z_2+i(c_2z_1-c_1z_2)^3+3+k)}{3}\right)\eeas is a solution of \eqref{e1.3a}.
\end{exm}
\begin{exm}
Let $a_1=1$, $a_2=-2$, $a_3=\sqrt{3}$ and $L(z)=z_1+2z_2$. Choose $k\in\mathbb{C}$ such that $e^{ik}=1$. Also, let $c=(c_1,c_2)\in\mathbb{C}^2$ such that $e^{iL(c)}=(1+i\sqrt{3})/2$. Then, it can be easily verified that \beas(f_1(z),f_2(z))=(\cos(z_1+2z_2+3),\cos(z_1+2z_2+3+k))\eeas is a solution of \eqref{e1.3a}.
\end{exm}
\begin{exm}
Let $a_1=1$, $a_2=-2$, $a_3=\sqrt{3}$ and $L(z)=5z_1-2z_2$. Choose $k\in\mathbb{C}$ such that $e^{ik}=-1$. Also, let $c=(c_1,c_2)\in\mathbb{C}^2$ such that $e^{iL(c)}=(1-i\sqrt{3})/2$. Then, it can be easily verified that \beas(f_1(z),f_2(z))=(\cos(5z_1-2z_2+10i),\cos(5z_1-2z_2+10i+k))\eeas is a solution of \eqref{e1.3a}.
\end{exm}
\begin{exm}
Let $a_1=12$, $a_2=7$, $a_3=5$ and $L(z)=z_1+iz_2$. Choose $k\in\mathbb{C}$ such that $e^{ik}=i$. Also, let $c=(c_1,c_2)\in\mathbb{C}^2$ such that $e^{iL(c)}=1$. Then, it can be easily verified that \beas(f_1(z),f_2(z))=\left(\frac{1}{12}\cos(-(z_1+iz_2+17i)+k),\frac{1}{12}\cos(z_1+iz_2+17i)\right)\eeas is a solution of \eqref{e1.3a}.
\end{exm}
\begin{theo}\label{t3}
Let $a_1,a_2,a_3, a_4$ be four non-zero constants in $\mathbb{C}$ such that $a_1^2+a_3^2=a_2+a_4^2$ and $D:=a_1a_4-a_2a_3\neq0$. Let $c=(c_1,c_2)\in \mathbb{C}^2\setminus\{0\}$. If $(f_1,f_2)$ is a pair of transcendental entire solution with finite order of Fermat-type simultaneous difference equation \eqref{e1.3b}, then $f_1(z)$ and $f_2(z)$ will be of the following form
\beas f_1(z)=\frac{1}{-2D}\left((a_3+ia_1e^{k})e^{{L(z)+A+\Phi(t)}}+(a_3-ia_1e^{-k})e^{-(L(z)+A+\Phi(t))}\right)\eeas and 
\beas f_2(z)=\frac{1}{-2D}\left((a_3e^{k}+ia_1)e^{{L(z)+A+\Phi(t)}}+(a_3e^{-k}-ia_1)e^{-(L(z)+A+\Phi(t))}\right),\eeas where $L(z)=\alpha_1z_1+\alpha_2z_2$, $e^{2k}=1$ with \beas e^{L(c)}=\frac{a_1+ia_3e^{k}}{a_2+ia_4e^{k}}=\frac{a_2e^{k}-ia_4}{a_1e^{k}-ia_3}=\frac{a_1+ia_3e^{-k}}{a_2+ia_4e^{-k}}=\frac{a_2-ia_4e^{k}}{a_1-ia_3e^{k}},\eeas $\alpha_1,\alpha_2, A,k\in\mathbb{C}$, and $\Phi(t)$ ia a polynomial in $t:=c_2z_1-c_1z_2$.
\end{theo}
\par The following examples show the existence of transcendental entire solutions with finite order of the system \eqref{e1.3b}.
\begin{exm}
Let us choose $a_1=a_2=a_3=1$, $a_4=-1$, $L(z)=z_1+2z_2$ and $c=(c_1,c_2)\in\mathbb{C}^2$ such that $c_1+2c_2=(2m+1/2)\pi$, $m$ being an integer. Let \beas f_1(z)=\frac{1}{4}\left((1+i)e^{z_1+2z_2+(c_2z_1-c_1z_2)^n+3}+(1-i)e^{-(z_1+2z_2+(c_2z_1-c_1z_2)^n+3)}\right)\eeas and \beas f_2(z)=\frac{1}{4}\left((1+i)e^{z_1+2z_2+(c_2z_1-c_1z_2)^n+3}-(i-1)e^{-(z_1+2z_2+(c_2z_1-c_1z_2)^n+3)}\right).\eeas Then one can easily verify that $(f_1(z),f_2(z))$ is a solution of \eqref{e1.3b}
\end{exm}
\begin{exm}
	Let us choose $a_1=a_2=a_3=1$, $a_4=-1$, $L(z)=i(z_1-z_2)$, $\Phi(t)=i(c_2z_1-c_1z_2)^5$, $A=3$ and $c=(c_1,c_2)\in\mathbb{C}^2$ such that $c_1+2c_2=(2m-1/2)\pi$, $m$ being an integer. Let \beas f_1(z)=\frac{\cos\left(z_1-z_2+(c_2z_1-c_1z_2)^5-3i\right)+\sin\left(z_1-z_2+(c_2z_1-c_1z_2)^5-3i\right)}{2}\eeas and \beas f_2(z)=-\frac{\cos\left(z_1-z_2+(c_2z_1-c_1z_2)^5-3i\right)+\sin\left(z_1-z_2+(c_2z_1-c_1z_2)^5-3i\right)}{2}\eeas Then one can easily verify that $(f_1(z),f_2(z))$ is a solution of \eqref{e1.3b}
\end{exm}
\par Besides finding the solutions of Fermat-type difference equations, Fermat-type partial differential-difference equations are also studied by many researchers (see \cite{Li-Zhang-Xu & 2021 & AIMS,Liu-Xu & JM & 2021,Xu & Cao & 2018,Xu & wang & 2020,Zheng-Xu & Analysis math & 2021}). For example, Xu and Cao \cite{Xu & Cao & 2018} have  investigated the entire solutions of Fermat-type partial differential-difference equation
\bea\label{e1.3} \left(\frac{\partial f(z_1,z_2)}{\partial z_1}\right)^n+f^m(z_1+c_1,z_2+c_2)=1\eea and obtained the following interesting result for the functions in $ \mathbb{C}^2 $.

\begin{theoG}\cite{Xu & Cao & 2018}
Let $c=(c_1,c_2)$ be a constant in $\mathbb{C}^2$. Then the Fermat-type partial differential-difference equation $(\ref{e1.3})$ does not have any transcendental entire solution with finite order, where $m$ and $n$ are two distinct positive integers.
\end{theoG}
In $2020$, Xu and Wang \cite{Xu & wang & 2020} generalized Theorem G by considering the following Fermat-type partial differential-difference equation \bea\label{e1.4} \left(\frac{\partial f(z_1,z_2)}{\partial z_1}+\frac{\partial f(z_1,z_2)}{\partial z_2}\right)^n+f^m(z_1+c_1,z_2+c_2)=1\eea
and proved the following result.
\begin{theoH}\cite{Xu & wang & 2020}
Let $c=(c_1,c_2)$ be a constant in $\mathbb{C}^2$ and $ m $, $ n $ be two positive integers. If the Fermat-type partial differential-difference equation \eqref{e1.4} satisfies one of the following conditions:
\begin{enumerate}
\item[(i)] $ m>n $;
\item[(ii)] $ n>m\geq 2 $, 
\end{enumerate}
then \eqref{e1.4} does not have any finite order transcendental entire solutions.
\end{theoH} 
Corresponding to Theorem G, Xu \textit{et. al.} \cite{Xu-Liu-Li-JMAA-2020} considered system of partial differential-difference equations and obtained the result as follows.
\begin{theoI}\cite{Xu-Liu-Li-JMAA-2020}
Let $c=(c_1,c_2)$ be a constant in $\mathbb{C}^2$, and $m_j,n_j$ $(j=1,2)$ be positive integers. If the following system of Fermat-type partial differential-difference equations
\beas \begin{cases} \left(\frac{\partial f_1(z_1,z_2)}{\partial z_1}\right)^{n_1}+f_2(z_1+c_1,z_2+c_2)^{m_1}=1\\\left(\frac{\partial f_2(z_1,z_2)}{\partial z_1}\right)^{n_2}+f_1(z_1+c_1,z_2+c_2)^{m_2}=1,\end{cases}\eeas satisfies one of the conditions
\begin{enumerate}
\item [(i)] $m_1m_2>n_1n_2$; \\
\item [(ii)] $m_j>\displaystyle\frac{n_j}{n_j-1}$, $j=1,2$,\\ then the above system does not have any pair of transcendental entire solution with finite order.\end{enumerate}
\end{theoI}
\par Liu and Xu \cite{Liu-Xu & JM & 2021}, further extended Theorem I by considering Fermat-type systems of second-order partial differential-difference equation, obtained the following result.
\begin{theoJ}\cite{Liu-Xu & JM & 2021}
Let $c=(c_1,c_2)\in\mathbb{C}^2$, and $m_j,n_j$ $(j=1,2)$ be positive integers, and $\alpha, \beta$ be constants in $\mathbb{C}$ that are not zero at the same time. If the following system of Fermat-type
partial differential-difference equations
\beas \begin{cases} \left(\alpha\frac{\partial^2 f_1(z_1,z_2)}{\partial z_1^2}+\beta\frac{\partial^2 f_1(z_1,z_2)}{\partial z_2^2}\right)^{n_1}+f_2(z_1+c_1,z_2+c_2)^{m_1}=1\\\left(\alpha\frac{\partial^2 f_2(z_1,z_2)}{\partial z_1^2}+\beta\frac{\partial^2 f_2(z_1,z_2)}{\partial z_2^2}\right)^{n_2}+f_1(z_1+c_1,z_2+c_2)^{m_2}=1,\end{cases}\eeas satisfies one of the conditions
\begin{enumerate}
\item [(i)] $m_1m_2>n_1n_2$; \\
\item [(ii)] $m_j>\displaystyle\frac{n_j}{n_j-1}$, $j=1,2$,\\ then the above system does not have any pair of transcendental entire solution with finite order.\end{enumerate}
\end{theoJ}
\par As far as we know, it appears that the Fermat-type mixed partial differential-difference equations in several complex variables has not been addressed in the literature before. In order to generalize and also to establish a result which combines Theorem I and Theorem J, we consider the following partial differential-difference equation
\bea\label{e22.5}\begin{cases} (a\partial^{I}f_1(z_1,z_2)+b \partial^{J}f_1(z_1,z_2))^{n_1}+f_2(z_1+c_1,z_2+c_2)^{m_1}=1\\(a \partial^{I}f_2(z_1,z_2)+b \partial^{J}f_2(z_1,z_2))^{n_2}+f_1(z_1+c_1,z_2+c_2)^{m_2}=1,\end{cases}\eea where 
\beas \partial^{I}f_j(z_1,z_2)=\displaystyle\frac{\partial^{|I|}f_j(z_1,z_2)}{\partial z_1^{\alpha_1}\partial z_2^{\alpha_2}}\;\; \mbox{and}\;\; \partial^{J}f_j(z_1,z_2)=\displaystyle\frac{\partial^{|J|}f_j(z_1,z_2)}{\partial z_1^{\beta_1}\partial z_2^{\beta_2}} \eeas with $I=(\alpha_1,\alpha_2)$ and $J=(\beta_1,\beta_2)$ are two multi-index, where $\alpha_1,\alpha_2$, $\beta_1$ and $\beta_2$ are non-negative integers and $a,b\in\mathbb{C}$, not both zero. We denote by $\mid I\mid$ to denote the length of $I$, that is, $\mid I\mid=\alpha_1+\alpha_2$. Similarly, for $J$ also. \vspace{1.2mm}
\par As a matter of fact, we prove the next result for any order Fermat-type partial differential-difference equation \eqref{e22.5}.
\begin{theo}\label{t4}
Let $c=(c_1,c_2)$ be a constant in $\mathbb{C}^2$ and $m_j$, $n_j $ be positive integers with $j=1,2$. If the Fermat-type simultaneous partial differential-difference equation \eqref{e22.5} satisfies one of the following conditions:
\begin{enumerate}
\item[(i)] $ m_1m_2>n_1n_2$;
\item[(ii)] $ m_j>\displaystyle\frac{n_j}{n_j-1}$, for $n_j\geq2$, $j=1,2$ 
\end{enumerate}
then \eqref{e22.5} does not have any pair of finite order transcendental entire solutions of the form $(f_1,f_2)$.
\end{theo}
\section{\textbf{Key Lemmas}} 
In this section, we present some necessary lemmas which will play key role to prove the main results of this paper.
\begin{lem}\label{lem3.1}\cite{Hu & Li & Yang & 2003}
Let $f_j\not\equiv0$ $(j=1,2\ldots,m;\; m\geq 3)$ be meromorphic functions on $\mathbb{C}^{n}$ such that $f_1,\ldots, f_{m-1}$ are not constants, $f_1+f_2+\cdots+ f_m=1$ and such that
\beas \sum_{j=1}^{m}\left\{N_{n-1}\left(r,\frac{1}{f_j}\right)+(m-1)\ol N(r,f_j)\right\}< \lambda T(r,f_j)+O(\log^{+}T(r,f_j))\eeas
holds for $j=1,\ldots, m-1$ and all $r$ outside possibly a set with finite logarithmic measure, where $\lambda < 1$ is a positive number. Then $f_m = 1$.	
\end{lem}
\begin{lem}\label{lem3.2}\cite{Lelong & 1968, ronkin & AMS & 1971, Stoll & AMS & 1974}
For an entire function $F$ on $\mathbb{C}^n$, $F(0)\not\equiv 0$ and put $\rho(n_F)=\rho<\infty$. Then there exist a canonical function $f_F$ and a function $g_F\in\mathbb{C}^n$ such that $F(z)=f_F (z)e^{g_F(z)}$. For the special case $n=1$, $f_F$ is the canonical product of Weierstrass.
\end{lem}
\begin{lem}\label{lem3.3}\cite{P`olya & Lond & 1926}
If $g$ and $h$ are entire functions on the complex plane $\mathbb{C}$ and $g(h)$ is an entire function of finite order, then there are only two possible
cases: either
\begin{enumerate}
\item [(i)] the internal function $h$ is a polynomial and the external function $g$ is of finite order; or else
\item [(ii)] the internal function $h$ is not a polynomial but a function of finite order, and the external function $g$ is of zero order.
\end{enumerate}\end{lem}

\begin{lem}\label{lem3.4}\cite{Biancofiore & Stoll & 1981, Ye & 1995}
Let $f$ be a non-constant meromorphic function on $\mathbb{C}^n$ and let $I=(\alpha_1,\ldots,\alpha_n)$ be a multi-index with length $|I|=\sum_{j=1}^{n}\alpha_j$. Assume that
$T(r_0,f)\geq e$ for some $r_0$. Then
\beas m\left(r,\frac{\partial^{I}f}{f}\right)=S(r,f)\eeas holds for all $r\geq r_0$ outside a set $E\subset
(0,+\infty) $ of finite logarithmic measure, $\displaystyle\int_{E}\frac{dt}{t}< \infty$, where $\partial^{I}f=\displaystyle\frac{\partial^{|I|}f}{\partial z_1^{\alpha_1}\ldots\partial z_2^{\alpha_2}}$.  
\end{lem}
\begin{lem}\label{lem2.5}\cite{Cao & Korhonen & 2016, Korhonen & CMFT & 2012}
Let $f$ be a non-constant meromorphic function with finite order on $\mathbb{C}^n$ such that $f(0)\neq 0,\infty$, and let $\epsilon>0$. Then for $c\in \mathbb{C}^n$,  \beas m\left(r,\frac{f(z+c)}{f(z)}\right)+m\left(r,\frac{f(z)}{f(z+c)}\right)=S(r,f)\eeas holds for all $r\geq r_0$ outside a set $E\subset
(0,+\infty) $ of finite logarithmic measure, $\displaystyle\int_{E}\frac{dt}{t}< \infty$.\end{lem}
\section{Proof of the main results}

\begin{proof}[\textbf{Proof of Theorem $\ref{t1}$}]
Suppose that $(f_1,f_2)$ is a pair of transcendental entire functions with finite order satisfying system \eqref{e1.2}. We write \eqref{e1.2} as the following:

\bea\label{e3.1} \begin{cases}
(f_1(z)+i(f_2(z+c)-f_2(z)))(f_1(z)-i(f_2(z+c)-f_2(z)))=1\\(f_2(z)+i(f_1(z+c)-f_1(z)))(f_2(z)-i(f_1(z+c)-f_1(z)))=1.
\end{cases}\eea
\par From the above equations we see that $f_i(z)\pm i(f_j(z+c)-f_j(z))$ have no zeros in $\mathbb{C}^2$, where $i\neq j$, $i,j=1,2$.\vspace{1mm}

\par Since $f_1$, $f_2$ are transcendental entire functions with finite order, there exist polynomials $p_1(z)$, $p_2(z)$ such that
\bea\label{e3.2} \begin{cases}
f_1(z)+ i(f_2(z+c)-f_2(z))=e^{p_1(z)}\\f_1(z)- i(f_2(z+c)-f_2(z))=e^{-p_1(z)}\\f_2(z)+ i(f_1(z+c)-f_1(z))=e^{p_2(z)}\\f_2(z)- i(f_1(z+c)-f_1(z))=e^{-p_2(z)}.\end{cases}\eea

\par In view of \eqref{e3.2}, we obtain 
\bea\label{e3.3} \begin{cases}
	f_1(z)=\displaystyle\frac{1}{2}\left(e^{p_1(z)}+e^{-p_1(z)}\right)\\ f_2(z+c)-f_2(z)=\displaystyle\frac{1}{2i}\left(e^{p_1(z)}-e^{-p_1(z)}\right)\\f_2(z)=\displaystyle\frac{1}{2}\left(e^{p_2(z)}+e^{-p_2(z)}\right)\\f_1(z+c)-f_1(z)=\displaystyle\frac{1}{2i}\left(e^{p_2(z)}-e^{-p_2(z)}\right).\end{cases}\eea
\par Form \eqref{e3.3}, we can easily obtain the following two equations:
\bea\label{e3.4} -ie^{p_2(z)+p_1(z+c)}-ie^{p_2(z)-p_1(z+c)}+ie^{p_2(z)+p_1(z)}+ie^{p_2(z)-p_1(z)}+e^{2p_2(z)}=1\eea and 

\bea\label{e3.5} -ie^{p_1(z)+p_2(z+c)}-ie^{p_1(z)-p_2(z+c)}+ie^{p_1(z)+p_2(z)}+ie^{p_1(z)-p_2(z)}+e^{2p_1(z)}=1.\eea

Now we consider the following two cases.\vspace{1mm}
\par \textbf{Case 1:} Suppose $p_2(z)-p_1(z)$ is constant. Let $p_2(z)-p_1(z)=k$, where $k\in\mathbb{C}$.\vspace{1mm}
\par Then \eqref{e3.4} and \eqref{e3.5}, respectively yield 
\bea\label{e3.6} -ie^{p_1(z)+p_1(z+c)+k}-ie^{p_1(z)-p_1(z+c)+k}+ie^{2p_1(z)+k}+e^{2p_1(z)+2k}=1-ie^{k}\eea 

and \bea\label{e3.7} -ie^{p_1(z)+p_1(z+c)+k}-ie^{p_1(z)-p_1(z+c)-k}+ie^{2p_1(z)+k}+e^{2p_1(z)}=1-ie^{-k}.\eea
\par First we consider that $e^{k}\neq\pm i$. Then, after simple calculation, we obtain from \eqref{e3.6} and \eqref{e3.7} that 
\bea\label{e3.8} \frac{-ie^k}{1-ie^k}e^{p_1(z)+p_1(z+c)}+\frac{-ie^k}{1-ie^k}e^{p_1(z)-p_1(z+c)}+\frac{(i+e^k)e^k}{1-ie^k}e^{2p_1(z)}=1\eea and 
\bea\label{e3.9} \frac{-ie^k}{1-ie^{-k}}e^{p_1(z)+p_1(z+c)}+\frac{-ie^{-k}}{1-ie^{-k}}e^{p_1(z)-p_1(z+c)}+\frac{(1+ie^k)e^k}{1-ie^{-k}}e^{2p_1(z)}=1.\eea
Using Lemma \ref{lem3.1}, we obtain from \eqref{e3.8} that \bea\label{e3.10} \frac{-ie^{k}}{1-ie^{k}}e^{p_1(z)-p_1(z+c)}=1.\eea

\par From \eqref{e3.8} and \eqref{e3.10}, we get
\bea\label{e3.11} e^{-p_1(z)+p_1(z+c)}=1-ie^{k}.\eea
\par Similarly, using Lemma \ref{lem3.1}, we obtain from \eqref{e3.9} that 
\bea\label{e3.12} \frac{-ie^{-k}}{1-ie^{-k}}e^{p_1(z)-p_1(z+c)}=1.\eea
Multiplying \eqref{e3.11} and \eqref{e3.12}, we get $-ie^{-k}=(1-ie^{-k})(1-ie^k)$, which yields that $ie^{k}=0$, which is not possible.\vspace{1mm}

\par Next, suppose that $e^{k}=-i$. Then from \eqref{e3.6}, we obtain that 
\beas e^{p_1(z)+p_1(z+c)}+e^{p_1(z)-p_1(z+c)}-2e^{2p_1(z)}=0,\eeas
which in turns into \beas e^{-p_1(z)+p_1(z+c)}+e^{-(p_1(z)+p_1(z+c))}=-2,\eeas which is not possible since the L.H.S. is transcendental entire, whereas the R.H.S. is constant.\vspace{1mm}
\par If $e^{k}=i$, then after simplification, we obtain from \eqref{e3.7} that
\beas e^{-p_1(z)+p_1(z+c)}-e^{-(p_1(z)+p_1(z+c))}=-2,\eeas which is again a contradiction by the same reason as discussed for the case $e^{k}=-i$.\vspace{1mm}
\par \textbf{Case 2:} Suppose $p_2(z)-p_1(z)$ is non-constant.\vspace{1mm}
\par We claim that $p_2(z)-p_1(z+c)$ is non-constant. If not, then $p_2(z)-p_1(z+c)=\textit{constant}=k$, say, where $k$ is a complex constant in one variable.\vspace{1mm}
\par This implies that $p_2(z)=L(z)+A$, where $L(z)=a_1z_1+a_2z_2$, $a_1,a_2, A\in \mathbb{C}$, and hence $p_1(z+c)=L(z)+A-k$. From this we get $p_1(z)=L(z)+A-L(c)-k$.\vspace{1mm}
\par Therefore, we must have $p_2(z)-p_1(z)=L(c)+k=\textit{constant}$, which is a contradiction.\vspace{1mm}
\par Now we consider the following subcases:\vspace{1mm}
\par \textbf{Subcase 2.1:} Let $p_2(z)+p_1(z)$ be constant, say $k$, where $k\in\mathbb{C}$.\vspace{1mm}
 \par Then, it is easily seen that both $p_2(z)-p_1(z+c)$ and $p_1(z)-p_2(z+c)$ are non-constants. \vspace{1mm}
 
\par Therefore, from \eqref{e3.3} and \eqref{e3.5}, we obtain 
\bea\label{e3.13} -ie^{p_2(z)+p_1(z+c)}-ie^{p_2(z)-p_1(z+c)}+ie^{p_2(z)-p_1(z)}+e^{2p_2(z)}=1-ie^k\eea and 
\bea\label{e3.14} -ie^{p_1(z)+p_2(z+c)}-ie^{p_1(z)-p_2(z+c)}+ie^{p_1(z)-p_2(z)}+e^{2p_1(z)}=1-ie^k,\eea respectively.\vspace{1mm}
\par Suppose $e^{k}\neq -i$. Then \eqref{e3.13} and \eqref{e3.14} can be written as
\bea\label{e3.15}  \left(\frac{-i}{1-ie^k}\right)e^{p_2(z)+p_1(z+c)}&+&\left(\frac{-i}{1-ie^k}\right)e^{p_2(z)-p_1(z+c)}+\left(\frac{i}{1-ie^{k}}\right)e^{p_2(z)-p_1(z)}\nonumber\\&+&\left(\frac{1}{1-ie^{k}}\right)e^{2p_2(z)}=1\eea and 
\bea\label{e3.16}  \left(\frac{-i}{1-ie^k}\right)e^{p_1(z)+p_2(z+c)}&+&\left(\frac{-i}{1-ie^k}\right)e^{p_1(z)-p_2(z+c)}+\left(\frac{i}{1-ie^{k}}\right)e^{p_1(z)-p_2(z)}\nonumber\\&+&\left(\frac{1}{1-ie^{k}}\right)e^{2p_1(z)}=1.\eea
\par Using Lemma \ref{lem3.1}, we obtain from \eqref{e3.15} that 
\bea\label{e3.17} \left(\frac{-i}{1-ie^k}\right)e^{p_2(z)+p_1(z+c)}=1.\eea
\par From \eqref{e3.15} and \eqref{e3.17}, we get after simple calculation that 
\bea\label{e3.18} e^{-(p_2(z)+p_1(z+c))}=-i(1+ie^{-k}).\eea
\par Using Lemma \ref{lem3.1}, we obtain from \eqref{e3.15} that 
\bea\label{e3.19} \left(\frac{-i}{1-ie^k}\right)e^{p_1(z)+p_2(z+c)}=1.\eea
\par From \eqref{e3.16} and \eqref{e3.19}, we get after simple calculation that 
\bea\label{e3.20} e^{-(p_1(z)+p_2(z+c))}=-i(1+ie^{-k}).\eea
\par From \eqref{e3.17} and \eqref{e3.18}, we obtain after simple calculation
\bea\label{e3.21} ie^{2k}-e^k-i=0.\eea
\par Observing \eqref{e3.17}, we conclude that $p_2(z)+p_1(z+c)$ must be constant, and hence we may assume that $p_2(z)=L(z)+A$, where $L(z)=a_1z_1+a_2z_2$, $a_1,a_2, A\in\mathbb{C}$.\vspace{1mm}
\par From \eqref{e3.17} and \eqref{e3.19}, we obtain 
\bea\label{e3.22} \left(\frac{-i}{1-ie^k}\right)e^{-L(z)+k}=1\eea and 
\bea\label{e3.23} \left(\frac{-i}{1-ie^k}\right)e^{L(z)+k}=1.\eea
\par Multiplying \eqref{e3.22} and \eqref{e3.23}, we get $e^{k}=1/2i$, which does not satisfy the relation \eqref{e3.21}.\vspace{1mm}
\par If $e^{k}=-i$, then after simple calculation, \eqref{e3.13} reduces to $e^{p_1(z)+p_1(z+c)}+e^{p_1(z)-p_1(z+c)}=-2$, which is again a contradiction since, the L.H.S. of the equation is transcendental entire, whereas the R.H.S. of it is constant.\vspace{1mm}
\par \textbf{Subcase 2.2:} Suppose $p_2(z)+p_1(z)$ is non-constant.\vspace{1mm}
\par We claim that $p_2(z)+p_1(z+c)$ is constant. If not, then by Lemma \ref{lem3.1} and \eqref{e3.4}, we obtain that $-ie^{p_2(z)-p_1(z+c)}=1$, which yields $p_2(z)-p_1(z+c)$ is constant, say $k$, where $k\in\mathbb{C}$.\vspace{1mm}
\par Thus, we may assume that $p_2(z)=L(z)+A,$ where $L(z)=a_1z_1+a_2z_2$, $a_1,a_2,A\in\mathbb{C}$. Therefore, $p_1(z+c)=L(z)+A-k$ and hence $p_1(z)=L(z)+A-L(c)-k$.\vspace{1mm} 
\par This implies that $p_2(z)-p_1(z)=L(c)+k$, a constant, which contradicts to our assumption. \vspace{1mm}
\par Let $p_2(z)+p_1(z+c)=k_1$, a complex constant in one variable.\vspace{1mm}

\par Then, \eqref{e3.4} reduces to 
\bea\label{e3.24} -ie^{p_2(z)-p_1(z+c)}+ie^{p_2(z)+p_1(z)}+ie^{p_2(z)-p_1(z)}+e^{2p_2(z)}=1+ie^{k_1}.\eea
\par Suppose $e^{k_1}\neq i$. Then using Lemma \ref{lem3.1}, we obtain from \eqref{e3.24} that \beas\left(\frac{-i}{1+ie^{k_1}}\right)e^{p_2(z)-p_1(z+c)}=1.\eeas

\par But this implies that $e^{p_2(z)-p_1(z+c)}$ and hence $p_2(z)-p_1(z+c)$ must be constant, say $k_2$, where $k_2\in\mathbb{C}$.\vspace{1mm}
\par Therefore, we must have $p_2(z)=(k_1+k_2)/2$, a constant, which contradicts to the fact that $p_2(z)$ is non-constant polynomial.\vspace{1mm}
\par If $e^{k_1}=i$, then after simplification, \eqref{e3.24} reduces to $e^{2p_1(z)}=-1$. But, then $p_1(z)$ becomes constant, which is a contradiction.\vspace{1mm}
\par \textbf{Subcase 2.3:} Suppose $p_2(z)+p_1(z+c)$ is a constant, say $k$, where $k\in \mathbb{C}$.\vspace{1mm} 
\par Then we may assume $p_2(z)=L(z)+A$, where $L(z)=a_1z_1+a_2z_2$, $a_1,a_2,A$ are complex constants in one variable, and hence $p_1(z+c)=-(L(z)+A)+k$.\vspace{1mm} 
\par Then,  it is clear that $p_1(z)=-(L(z)+A)+L(c)+k$, and therefore, $p_2(z)+p_1(z)=L(c)+k$ and $p_1(z)+p_2(z+c)=2L(c)+k$.\vspace{1mm}
With these, \eqref{e3.4} and \eqref{e3.5}, respectively reduce to 
\bea\label{e3.25} -ie^{p_2(z)-p_1(z+c)}+ie^{p_2(z)-p_1(z)}+e^{2p_2(z)}=1+ie^{k}\left(1-e^{L(c)}\right)\eea and
\bea\label{e3.26} -ie^{p_1(z)-p_2(z+c)}+ie^{p_1(z)-p_2(z)}+e^{2p_1(z)}=1- ie^{L(c)+k}\left(1-e^{L(c)}\right).\eea
\par First suppose $1+ie^{k}\left(1-e^{L(c)}\right)\neq0$.\vspace{1mm}
\par Then, using Lemma \ref{lem3.1}, we obtain from \eqref{e3.25} that 
\beas \left(\frac{-i}{1+ie^{k}(1-e^{L(c)})}\right)e^{p_2(z)-p_1(z+c)}=1,\eeas 
which implies that $p_2(z)-p_1(z+c)$ is constant. But This is not possible since $p_2(z)-p_1(z+c)=2(L(z)+A)-k$.\vspace{1mm}
\par Next Suppose that $1+ie^{k}\left(1-e^{L(c)}\right)=0$.\vspace{1mm}
\par Then after simple calculation, \eqref{e3.25} yields \beas ie^{-(p_2(z)+p_1(z+c))}-ie^{-(p_2(z)+p_1(z))}=1,\eeas
which implies that 
\bea\label{e3.27} ie^{-k}\left(1-e^{-L(c)}\right)=1.\eea
Also,\eqref{e3.26} yields 
\bea\label{e3.28} -ie^{p_1(z)-p_2(z+c)}+ie^{p_1(z)-p_2(z)}+e^{2p_1(z)}=1+e^{L(c)}.\eea 
Now, if $1+e^{L(c)}\neq0$, then using Lemma \ref{lem3.1} to \eqref{e3.28}, we obtain \beas \left(\frac{-i}{1+e^{L(c)}}\right)e^{p_1(z)-p_2(z+c)}=1,\eeas which implies that $p_1(z)-p_2(z+c)$ is constant, which is not possible since $p_1(z)-p_2(z+c)=-2(L(z)+A)+k$.\vspace{1mm}
\par So, $e^{L(c)}=-1$. Therefore, putting the value of $e^{L(c)}$ to $1+ie^{k}\left(1-e^{L(c)}\right)=0$, we obtain $e^{-k}=2i$. Hence, from \eqref{e3.27}, we get $-4=1$, which is a contradiction.\vspace{1mm}

\par \textbf{Subcase 2.4:} Suppose $p_2(z)+p_1(z+c)$ is non-constant.\vspace{1mm}
\par If $p_2(z)-p_1(z+c)$ is constant, say $B$, where $B\in\mathbb{C}$, then we may assume $p_2(z)=L(z)+A$, where $L(z)=a_1z_1+a_2z_2$ with $a_1,a_2,A\in \mathbb{C}$. Therefore, we have $p_1(z+c)=L(z)+A-B$, and hence $p_1(z)=L(z)+A-L(c)-B$.\vspace{1mm}
\par This implies that $p_2(z)-p_1(z)=L(c)+B$, a constant, which is a contradiction. Hence, $p_2(z)-p_1(z+c)$ is non-constant.\vspace{1mm}
\par Therefore, using Lemma \ref{lem3.1}, we obtain from \eqref{e3.4} that 
\bea\label{e3.29} ie^{(p_2(z)+p_1(z))}=1.\eea
\par Using \eqref{e3.29} in \eqref{e3.4}, we get 
\bea\label{e3.30} ie^{(-p_2(z)+p_1(z+c))}+ie^{-(p_2(z)+p_1(z+c))}-ie^{-(p_2(z)+p_1(z))}=1.\eea

\par From \eqref{e3.29}, we conclude that $p_2(z)+p_1(z)=\textit{constant}=k$, say.\vspace{1mm}
\par Using this in \eqref{e3.30}, we get 
\beas e^{-p_2(z)+p_1(z+c)}+e^{-(p_2(z)+p_1(z+c))}=-i(1+ie^{-k}), \eeas which is not possible since, the L.H.S. is transcendental entire whereas the R.H.S. is constant.\vspace{1mm}
\par Hence, we conclude that there is no pair of transcendental entire solutions with finite order of the system \eqref{e1.2}.\vspace{1mm}
\end{proof}
\begin{proof}[\textbf{Proof of Theorem $\ref{t2}$}]
As we all know the fact that the entire solutions of the functional equation $f^2+g^2=1$ are $f=\cos \alpha(z)$ and $g=\sin \alpha(z)$, where $\alpha(z)$ is an entire function. If $f,\; g$ are finite order entire functions, then $p(z)$ must be a non-constant polynomial (see \cite{Gross & Bull. Amer. & 1966, Gross & Amer. Math. & 1966, Montel & Paris & 1927}).\vspace{1mm}
\par Keeping in view of the above fact, we may obtain from \eqref{e1.3a} that \bea\label{e3.31}  \begin{cases}
a_1f_1(z)=\displaystyle\frac{1}{2}\left(e^{ip_1(z)}+e^{-ip_1(z)}\right)\\ a_2f_2(z+c)+a_3f_2(z)=\displaystyle\frac{1}{2i}\left(e^{ip_1(z)}-e^{-ip_1(z)}\right)\\a_1f_2(z)=\displaystyle\frac{1}{2}\left(e^{ip_2(z)}+e^{-ip_2(z)}\right)\\a_2f_1(z+c)+a_3f_1(z)=\displaystyle\frac{1}{2i}\left(e^{ip_2(z)}-e^{-ip_2(z)}\right),\end{cases}\eea where $p_1(z),\;p_2(z)$ are two non-constant polynomials.

\par After some simple calculations, we obtain from \eqref{e3.31} that 
\bea\label{e3.32} -\frac{ia_2}{a_1}e^{i(p_2(z)+p_1(z+c))}&-&\frac{ia_2}{a_1}e^{i(p_2(z)-p_1(z+c))}-\frac{ia_3}{a_1}e^{i(p_2(z)+p_1(z))}\nonumber\\&-&\frac{ia_3}{a_1}e^{i(p_2(z)-p_1(z))}+e^{2ip_2(z)}=1\eea

and \bea\label{e3.33} -\frac{ia_2}{a_1}e^{i(p_1(z)+p_2(z+c))}&-&\frac{ia_2}{a_1}e^{i(p_1(z)-p_2(z+c))}-\frac{ia_3}{a_1}e^{i(p_1(z)+p_2(z))}\nonumber\\&-&\frac{ia_3}{a_1}e^{i(p_1(z)-p_2(z))}+e^{2ip_1(z)}=1.\eea

\par We now consider the following two cases.\vspace{1mm}

\par \textbf{Case 1:} Suppose $p_2(z)-p_1(z)$ is constant, say $k$, where $k\in\mathbb{C}$.\vspace{1mm}
\par Then, it can be easily seen that $p_2(z)+p_1(z)$ and $p_2(z)+p_1(z+c)$ are non-constant polynomials.\vspace{1mm}

\par Putting $p_2(z)-p_1(z)=k$ in \eqref{e3.32} and \eqref{e3.33}, we obtain

\bea\label{e3.34} -\frac{ia_2}{a_1}e^{i(p_2(z)+p_1(z+c))}-\frac{ia_2}{a_1}e^{i(p_2(z)-p_1(z+c))}&-&\frac{ia_3}{a_1}e^{i(p_2(z)+p_1(z))}+e^{2ip_2(z)}\nonumber\\&=&\frac{a_1+ia_3e^{ik}}{a_1}\eea

and \bea\label{e3.35} -\frac{ia_2}{a_1}e^{i(p_1(z)+p_2(z+c))}-\frac{ia_2}{a_1}e^{i(p_1(z)-p_2(z+c))}&-&\frac{ia_3}{a_1}e^{i(p_1(z)+p_2(z))}+e^{2ip_1(z)}\nonumber\\&=&\frac{a_1+ia_3e^{-ik}}{a_1},\eea respectively.\vspace{1mm}
\par If $e^{ik}\neq ia_1/a_3,\; -ia_3/a_1$, then \eqref{e3.34} and \eqref{e3.35}, respectively yield
 
\bea\label{e3.36}&& \left(\frac{-ia_2}{a_1+ia_3e^{ik}}\right)e^{i(p_2(z)+p_1(z+c))}+\left(\frac{-ia_2}{a_1+ia_3e^{ik}}\right)e^{i(p_2(z)-p_1(z+c))}\nonumber\\&+&\left(\frac{-ia_3}{a_1+ia_3e^{ik}}\right)e^{i(p_2(z)+p_1(z))}+\left(\frac{a_1}{a_1+ia_3e^{ik}}\right)e^{2ip_2(z)}=1\eea

and \bea\label{e3.37}&& \left(\frac{-ia_2}{a_1+ia_3e^{-ik}}\right)e^{i(p_1(z)+p_2(z+c))}+\left(\frac{-ia_2}{a_1+ia_3e^{-ik}}\right)e^{i(p_1(z)-p_2(z+c))}\nonumber\\&+&\left(\frac{-ia_3}{a_1+ia_3e^{-ik}}\right)e^{i(p_1(z)+p_2(z))}+\left(\frac{a_1}{a_1+ia_3e^{-ik}}\right)e^{2ip_1(z)}=1.\eea
\par Now, using Lemma \ref{e3.1}, we obtain from \eqref{e3.36} and \eqref{e3.37} that
\bea\label{e3.38} \left(\frac{-ia_2}{a_1+ia_3e^{ik}}\right)e^{i(p_2(z)-p_1(z+c))}=1\eea and \bea\label{e3.39} \left(\frac{-ia_2}{a_1+ia_3e^{-ik}}\right)e^{i(p_1(z)-p_2(z+c))}=1,\eea respectively.
\par Keeping in view of \eqref{e3.38}, \eqref{e3.36} reduces to \beas ia_2e^{i(p_2(z)+p_1(z+c))}+ia_3e^{i(p_2(z)+p_1(z))}=a_1e^{2ip_2(z)}.\eeas i.e., \beas ia_2e^{i(-p_2(z)+p_1(z+c))}+ia_3e^{i(-p_2(z)+p_1(z))}=a_1.\eeas i.e., \bea\label{e3.40} ia_2e^{i(-p_2(z)+p_1(z+c))}=a_1-ia_3e^{-ik}.\eea
\par Similarly, using \eqref{e3.39} in \eqref{e3.37}, we obtain
\bea\label{e3.41} ia_2e^{i(-p_1(z)+p_2(z+c))}=a_1-ia_3e^{ik}.\eea
\par Now, we divide \eqref{e3.38} by \eqref{e3.39} to get \beas e^{i(p_2(z)-p_1(z)+p_2(z+c)-p_1(z+c))}=1.\eeas
\par As $p_2(z)-p_1(z)=k$, we obtain from the above equation \bea\label{e3.42} e^{2ik}=1.\eea
\par Also, multiplying \eqref{e3.38} and \eqref{e3.40}, we obtain with the help of \eqref{e3.42} that
\beas a_2^2=(a_1+ia_3e^{ik})(a_1-ia_3e^{-ik}).\eeas i.e., \beas a_2^2=a_1^2+a_3^2-ia_1a_3e^{-ik}+ia_1a_3e^{ik}.\eeas i.e., \bea\label{e3.43} a_2^2=a_1^2+a_3^2.\eea
\par Since, $p_2(z)-p_1(z)=k$, \eqref{e3.38} reduces to \beas\left(\frac{-ia_2}{a_1+ia_3e^{ik}}\right)e^{i(p_1(z)-p_1(z+c)+k)}=1,\eeas which implies that $e^{i(p_1(z)-p_1(z+c)+k)}$ and hence $p_1(z)-p_1(z+c)$ must be constant. 
\par Therefore, we may assume that $p_1(z)=L(z)+\Phi(t)+A,$ where $L(z)=\alpha_1z_1+\alpha_2z_2$ with $\alpha_1,\;\alpha_2,\; A\in \mathbb{C}$ and $\Phi(t)$ is a polynomial in $t:=c_2z_1-c_1z_2$.\vspace{1mm}

\par Now, we divide \eqref{e3.38} by \eqref{e3.40} to get \beas e^{2i(p_2(z)-p_1(z+c))}=-\frac{a_1+ia_3e^{ik}}{a_1-ia_3e^{-ik}}.\eeas 
\par Since $p_2(z)-p_1(z)=k$, $p_1(z)=L(z)+\Phi(t)+A$ and $e^{2ik}=1$, we obtain from the above equation that \beas e^{2iL(c)}=-\frac{a_1-ia_3e^{-ik}}{a_1+ia_3e^{ik}}.\eeas
\par Therefore, the solution of the system \eqref{e1.3a} is 
\beas(f_1(z),f_2(z))=\displaystyle\left(\frac{1}{a_1}\cos(L(z)+\Phi(t)+A),\displaystyle\frac{1}{a_1}\cos(L(z)+\Phi(t)+A+k)\right).\eeas

\par \textbf{Case 2:} Suppose $p_2(z)-p_1(z)$ is non-constant.\vspace{1mm}
\par We claim that $p_2(z)-p_1(z+c)$ is non-constant. If not, suppose $p_2(z)-p_1(z+c)=k$, a constant in $\mathbb{C}$. Then we may assume that $p_2(z)=L(z)+\Phi(t)+A$, where $L(z)$, $\Phi(t)$, $A$ is defined as in Case 1. This implies that $p_1(z+c)=L(z)+\Phi(t)+A-k$ and hence $p_1(z)=L(z)+\Phi(t)+A-L(c)-k$. But, then $p_2(z)-p_1(z)=L(c)+k$, a constant, which contradicts to our assumption. So, our claim is true.\vspace{1mm} 
\par Now, we discuss the following subcases:\vspace{1mm}
\par \textbf{Subcase 2.1:} Let $p_2(z)+p_1(z)=k$, a constant, $k\in\mathbb{C}$. Then it can be easily verified that both $p_2(z)-p_1(z+c)$ and $p_1(z)-p_2(z+c)$ are non-constants.\vspace{1mm}

\par Therefore, from \eqref{e3.32} and \eqref{e3.33}, we obtain
\bea\label{e3.44} -\frac{ia_2}{a_1}e^{i(p_2(z)+p_1(z+c))}&-&\frac{ia_2}{a_1}e^{i(p_2(z)-p_1(z+c))}-\frac{ia_3}{a_1}e^{i(p_2(z)-p_1(z))}\nonumber\\&+&e^{2ip_2(z)}=\frac{1}{a_1}(a_1+ia_3e^{ik})\eea

and \bea\label{e3.45} -\frac{ia_2}{a_1}e^{i(p_1(z)+p_2(z+c))}&-&\frac{ia_2}{a_1}e^{i(p_1(z)-p_2(z+c))}-\frac{ia_3}{a_1}e^{i(p_1(z)-p_2(z))}\nonumber\\&+&e^{2ip_1(z)}=\frac{1}{a_1}(a_1+ia_3e^{ik}),\eea respectively.\vspace{1mm}

\par Suppose $a_1+ia_3e^{ik}\neq0$. Then \eqref{e3.44} and \eqref{e3.45}, respectively can be written as 
\bea\label{e3.46}&& \left(\frac{-ia_2}{a_1+ia_3e^{ik}}\right)e^{i(p_2(z)+p_1(z+c))}+\left(\frac{-ia_2}{a_1+ia_3e^{ik}}\right)e^{i(p_2(z)-p_1(z+c))}\nonumber\\&&+\left(\frac{-ia_3}{a_1+ia_3e^{ik}}\right)e^{i(p_2(z)-p_1(z))}+\left(\frac{a_1}{a_1+ia_3e^{ik}}\right)e^{2ip_2(z)}=1\eea and 
\bea\label{e3.47}&& \left(\frac{-ia_2}{a_1+ia_3e^{ik}}\right)e^{i(p_1(z)+p_2(z+c))}+\left(\frac{-ia_2}{a_1+ia_3e^{ik}}\right)e^{i(p_1(z)-p_2(z+c))}\nonumber\\&&+\left(\frac{-ia_3}{a_1+ia_3e^{ik}}\right)e^{i(p_1(z)-p_2(z))}+\left(\frac{a_1}{a_1+ia_3e^{ik}}\right)e^{2ip_1(z)}=1.\eea

\par Now, using Lemma \ref{lem3.1}, we obtain from \eqref{e3.46} that \bea\label{e3.48} \left(\frac{-ia_2}{a_1+ia_3e^{ik}}\right)e^{i(p_2(z)+p_1(z+c))}=1.\eea 
\par From \eqref{e3.46} and \eqref{e3.48}, we get
\bea\label{e3.49} ia_2e^{-i(p_2(z)+p_1(z+c))}=a_1-ia_3e^{-ik}.\eea
\par Again using Lemma \ref{lem3.1}, we obtain from \eqref{e3.47} that \bea\label{e3.50} \left(\frac{-ia_2}{a_1+ia_3e^{ik}}\right)e^{i(p_1(z)+p_2(z+c))}=1.\eea 
\par From \eqref{e3.46} and \eqref{e3.48}, we get
\bea\label{e3.51} ia_2e^{-i(p_1(z)+p_2(z+c))}=a_1-ia_3e^{-ik}.\eea

\par We observe that L.H.S. of \eqref{e3.48} is transcendental entire, whereas R.H.S. is constant.\vspace{1mm} 
\par Therefore, $e^{-i(p_2(z)+p_1(z+c))}$ and hence $p_2(z)+p_1(z+c)$ must be constant. Let $p_2(z)+p_1(z+c)=k_1$, where $k\in\mathbb{C}$. Then we may assume that $p_2(z)=L(z)+\Phi(t)+A$, where $L(z)=\alpha_1z_1+\alpha_2z_2$, $\alpha_1,\alpha_2, A\in \mathbb{C}$. This implies that $p_1(z+c)=-(L(z)+\Phi(t)+A)+L(c)+k_1$.\vspace{1mm}
\par Also, we may easily get that $p_1(z)=-(L(z)+\Phi(t)+A)+L(c)+k_1$, $p_2(z)+p_1(z)=L(c)+k_1=k$ and $p_1(z)+p_2(z+c)=2L(c)+k_1=L(c)+k$.\vspace{1mm}
\par Therefore, from \eqref{e3.31}, we obtain \beas (f_1(z),f_2(z))=\left(\frac{1}{a_1}\cos\left(-(L(z)+\Phi(t)+A)+k\right),\frac{1}{a_1}\cos(L(z)+\Phi(t)+A))\right).\eeas

\par Now, from \eqref{e3.49} and \eqref{e3.51}, we get $e^{ik_1}=e^{i(2L(c)+k_1)}$, which implies that \beas e^{2iL(c)}=1.\eeas 
\par Divide \eqref{e3.48} by \eqref{e3.49} to get \beas e^{2i(-L(c)+k)}=-\frac{a_1+ia_3e^{ik}}{a_1-ia_3e^{-ik}},\eeas
which implies that \beas e^{2ik}=-\frac{a_1+ia_3e^{ik}}{a_1-ia_3e^{-ik}}.\eeas\par After simplification, we obtain from the last equation that $e^{2ik}=-1$.\vspace{1mm}

\par If $e^{iL(c)}=1$ and $e^{ik}=i$, then from \eqref{e3.48}, we get $a_1=a_2+a_3$.\vspace{1mm}
\par If $e^{iL(c)}=1$ and $e^{ik}=-i$, then from \eqref{e3.48}, we get $a_1=-(a_2+a_3)$.\vspace{1mm}
\par If $e^{iL(c)}=-1$ and $e^{ik}=i$, then from \eqref{e3.48}, we get $a_1=-(a_2-a_3)$.\vspace{1mm}
\par If $e^{iL(c)}=-1$ and $e^{ik}=-i$, then from \eqref{e3.48}, we get $a_1=(a_2-a_3)$.\vspace{1mm}

\par If $a_1+ia_3e^{ik}=0$, then \eqref{e3.44} yields \beas -\frac{ia_2}{a_1}e^{i(p_2(z)+p_1(z+c))}-\frac{ia_2}{a_1}e^{i(p_2(z)-p_1(z+c))}-\frac{ia_3}{a_1}e^{i(p_2(z)-p_1(z))}\nonumber+e^{2ip_2(z)}=0.\eeas i.e.,
\beas ia_2e^{i(-p_2(z)+p_1(z+c))}+ia_2e^{-i(p_2(z)+p_1(z+c))}+ia_3e^{-i(p_2(z)+p_1(z))}=a_1.\eeas i.e., 
\beas ia_2\left(e^{i(-p_2(z)+p_1(z+c))}+e^{-i(p_2(z)+p_1(z+c))}\right)=a_1-ia_3e^{-ik},\eeas which implies that 
\beas e^{i(-p_2(z)+p_1(z+c))}+e^{-i(p_2(z)+p_1(z+c))}=\frac{-i(a_1^2-a_2^2)}{a_1a_2},\eeas which is not possible since L.H.S. is transcendental entire, whereas R.H.S. is constant.\vspace{1mm}

\par \textbf{Subcase 2.2:} Suppose $p_2(z)+p_1(z)$ is non-constant.\vspace{1mm}
\par We claim that $p_2(z)+p_1(z+c)$ is constant. If not, then by Lemma \ref{lem3.1}, we obtain from from \eqref{e3.32} that \beas -\frac{ia_2}{a_1}e^{i(p_2(z)-p_1(z+c))}=1,\eeas which implies that $p_2(z)-p_1(z+c)$ is constant, say $k\in\mathbb{C}$. Thus, we may assume that $p_2(z)=L(z)+\Phi(t)+A$, where $L(z)=\alpha_1z_1+\alpha_2z_2$, $\alpha_1,\alpha_2, A\in\mathbb{C}$ and $\Phi(t)$ is a polynomial in $t:=c_2z_1-c_1z_2$. This implies that $p_1(z+c)=(L(z)+\Phi(t)+A)-k$, and hence $p_1(z)=L(z)+\Phi(t)+A-k-L(c)$. But, then $p_2(z)-p_1(z)=L(c)+k=\text{constant}$, which is a contradiction. So, our claim is true.\vspace{1mm}
\par Suppose $p_2(z)+p_1(z+c)=k_1$, where $k_1\in\mathbb{C}$. Then \eqref{e3.32} reduces to 
\bea\label{e3.52} -\frac{ia_2}{a_1}e^{i(p_2(z)-p_1(z+c))}-\frac{ia_3}{a_1}e^{i(p_2(z)+p_1(z))}&-&\frac{ia_3}{a_1}e^{i(p_2(z)-p_1(z))}+e^{2ip_2(z)}\nonumber\\&=&\frac{1}{a_1}\left(a_1+ia_2e^{ik_1}\right).\eea
\par If $a_1+ia_2e^{ik_1}\neq0$, then using Lemma \ref{lem3.1}, we obtain from \eqref{e3.52} that 
\beas \left(\frac{-ia_2}{a_1+ia_3e^{ik_1}}\right)e^{i(p_2(z)-p_1(z+c))}=1.\eeas 
\par This implies that $e^{i(p_2(z)-p_1(z+c))}$ and hence $p_2(z)-p_1(z+c)$ is constant, say $k_2$, $k_2\in\mathbb{C}$. Therefore, we must have $p_2(z)=(k_1+k_2)/2=\text{constant}$, which contradicts to the fact that $p_2(z)$ is non-constant.\vspace{1mm}
\par If $a_1+ia_2e^{ik_1}=0$, then \eqref{e3.52} reduces to \beas -\frac{ia_2}{a_1}e^{i(p_2(z)-p_1(z+c))}-\frac{ia_3}{a_1}e^{i(p_2(z)+p_1(z))}&-&\frac{ia_3}{a_1}e^{i(p_2(z)-p_1(z))}+e^{2ip_2(z)}=0.\eeas i.e., 
\beas \frac{ia_2}{a_1}e^{-i(p_2(z)+p_1(z+c))}+\frac{ia_3}{a_1}e^{-i(p_2(z)-p_1(z))}&-&\frac{ia_3}{a_1}e^{-i(p_2(z)+p_1(z))}+e^{2ip_2(z)}=1.\eeas i.e.,
\beas ia_3\left(e^{-i(p_2(z)-p_1(z))}-e^{-i(p_2(z)+p_1(z))}\right)=a_1-ia_2e^{-ik_1}.\eeas i.e., 

\beas \left(e^{-i(p_2(z)-p_1(z))}-e^{-i(p_2(z)+p_1(z))}\right)=\frac{-i(a_1^2-a_2^2)}{a_1^2a_3},\eeas which is not possible since L.H.S. is transcendental entire, whereas R.H.S. is constant.\vspace{1mm}

\par \textbf{Subcase 2.3:} Suppose $p_2(z)+p_1(z+c)=\text{constant}=k$, say, $k\in\mathbb{C}$.\vspace{1mm}
\par Then, we may assume $p_2(z)=L(z)+\Phi(t)+A$, where $L(z)=\alpha_1z_1+\alpha_2z_2$, $\alpha_1,\;\alpha_2,\;A\in\mathbb{C}$ and $\Phi(t)$ is a polynomial in $t:=c_2z_1-c_1z_2$. So, $p_1(z+c)=-(L(z)+\Phi(t)+A)+k$, and hence $p_1(z)=-(L(z)+\Phi(t)+A)+L(c)+k$, $p_2(z)+p_1(z)=L(c)+k$ and $p_1(z)+p_2(z+c)=2L(c)+k$.\vspace{1mm}
\par Therefore, we have from \eqref{e3.31}
\beas (f_1(z),f_2(z))=\left(\frac{1}{a_1}\cos(-(L(z)+\Phi(t)+A)+L(c)+k),\frac{1}{a_1}\cos(L(z)+\Phi(t)+A)\right).\eeas
\par Now, \eqref{e3.32} reduces to
\bea\label{e3.53} -\frac{ia_2}{a_1}e^{i(p_2(z)-p_1(z+c))}&-&\frac{ia_3}{a_1}e^{i(p_2(z)-p_1(z))}+e^{2ip_2(z)}\nonumber\\&=&1+\frac{ia_2}{a_1}e^{ik}+\frac{ia_3}{a_1}e^{i(L(c)+k)}\nonumber\\&=&\frac{1}{a_1}\left(a_1+ie^{ik}\left(a_2+a_3e^{iL(c)}\right)\right).\eea 
\par Similarly, we obtain from \eqref{e3.33} that 
\bea\label{e3.54} -\frac{ia_2}{a_1}e^{i(p_1(z)-p_2(z+c))}&-&\frac{ia_3}{a_1}e^{i(p_1(z)-p_2(z))}+e^{2ip_1(z)}\nonumber\\&=&1+\frac{ia_2}{a_1}e^{ik}+\frac{ia_3}{a_1}e^{i(L(c)+k)}\nonumber\\&=&\frac{1}{a_1}\left(a_1+ie^{i(L(c)+k)}\left(a_2e^{iL(c)}+a_3)\right)\right).\eea 
\par If $a_1+ie^{ik}\left(a_2+a_3e^{iL(c)}\right)\neq0$, then using Lemma \ref{lem3.1}, we obtain from \eqref{e3.53}
\beas \frac{ia_2}{a_1+ie^{ik}(a_2+a_3e^{iL(c)})}e^{i(p_2(z)-p_1(z+c))}=1,\eeas which implies that $p_2(z)-p_1(z+c)$ is constant, a contradiction as it is non-constant.\vspace{1mm}
\par If $a_1+ie^{ik}\left(a_2+a_3e^{iL(c)}\right)=0$, then \eqref{e3.53} yields 
\beas -\frac{ia_2}{a_1}e^{i(p_2(z)-p_1(z+c))}&-&\frac{ia_3}{a_1}e^{i(p_2(z)-p_1(z))}+e^{2ip_2(z)}=0.\eeas i.e., 
\beas \frac{ia_2}{a_1}e^{-i(p_2(z)+p_1(z+c))}+\frac{ia_3}{a_1}e^{-i(p_2(z)+p_1(z))}=1.\eeas i.e., 
\bea\label{e3.55} ie^{-ik}\left(a_2+a_3e^{-iL(c)}\right)=a_1.\eea

\par Using \eqref{e3.55} in \eqref{e3.54}, we obtain
\bea\label{e3.56} -\frac{ia_2}{a_1}e^{i(p_1(z)-p_2(z+c))}&-&\frac{ia_3}{a_1}e^{i(p_1(z)-p_2(z))}+e^{2ip_1(z)}=1+e^{2i(L(c)+k)}.\eea 
\par If $1+e^{2i(L(c)+k)}\neq0$, then using Lemma \ref{lem3.1}, we obtain from that 
\beas \frac{-ia_2}{a_1(1+e^{2i(L(c)+k)})}e^{i(p_1(z)-p_2(z+c))}=1,\eeas which implies that $p_1(z)-p_2(z+c)$ is constant. But, $p_1(z)-p_2(z+c)=-2(L(z)+\Phi(t)+A)+2L(c)+k$, which is non-constant. So, we get a contradiction.\vspace{1mm}

\par So, it must be that \bea\label{e3.57}e^{2i(L(c)+k)}=-1.\eea 
\par Therefore, we obtain from \eqref{e3.56} that \bea\label{e3.58} ie^{-i(L(c)+k)}\left(a_3+a_2e^{-iL(c)}\right)=a_1.\eea
\par Now, comparing \eqref{e3.55} and \eqref{e3.58}, we get \beas e^{2iL(c)}=1.\eeas 
\par Hence, from \eqref{e3.57}, we get $e^{2ik}=-1.$\vspace{1mm}
\par If $e^{iL(c)}=1$ and $e^{ik}=i$, then from \eqref{e3.55}, we get $a_1=a_2+a_3$.\vspace{1mm}
\par If $e^{iL(c)}=1$ and $e^{ik}=-i$, then from \eqref{e3.55}, we get $a_1=-(a_2+a_3)$.\vspace{1mm}
\par If $e^{iL(c)}=-1$ and $e^{ik}=i$, then from \eqref{e3.55}, we get $a_1=a_2-a_3$.\vspace{1mm}
\par If $e^{iL(c)}=-1$ and $e^{ik}=-i$, then from \eqref{e3.55}, we get $a_1=-(a_2-a_3)$.\vspace{1mm}
\par \textbf{Subcase 2.4:} Suppose $p_2(z)+p_1(z+c)$ is non-constant.
Now, using Lemma \ref{lem3.1}, we obtain from \eqref{e3.32} that \beas -\frac{ia_3}{a_1}e^{i(p_1(z)+p_2(z))}=1,\eeas which implies that $p_2(z)+p_1(z)$ is constant, say $k$, where $k\in\mathbb{C}$.\vspace{1mm}
\par Using these in \eqref{e3.32}, we get \beas e^{-i(p_2(z)-p_1(z+c)+e^{-i(p_2(z)+p_1(z+c))})}=-\frac{i(a_1-ia_3e^{-ik})}{a_2},\eeas which is not possible since L.H.S. is transcendental entire, whereas R.H.S. is constant.\vspace{1mm}
\par Hence the proof.
\end{proof}
\begin{proof}[\textbf{Proof of Theorem $\ref{t3}$}]
Assume that $(f_1,f_2)$ is a pair of transcendental entire solution of \eqref{e1.3b} with each $f_j$ is of finite order, $j=1,2$.\vspace{1mm} 

\par Then, by similar argument as done in the proof of Theorem \ref{t2}, we obtain
\beas  \begin{cases}
a_1f_1(z+c)+a_2f_1(z)=\displaystyle\frac{1}{2}\left(e^{p_1(z)}+e^{-p_1(z)}\right)\vspace{1.5mm}\\ a_3f_2(z+c)+a_4f_2(z)=\displaystyle\frac{1}{2i}\left(e^{p_1(z)}-e^{-p_1(z)}\right)\vspace{1.5mm}\\a_1f_2(z+c)+a_2f_2(z)=\displaystyle\frac{1}{2}\left(e^{p_2(z)}+e^{-p_2(z)}\right)\vspace{1.5mm}\\a_3f_1(z+c)+a_4f_1(z)=\displaystyle\frac{1}{2i}\left(e^{p_2(z)}-e^{-p_2(z)}\right),\end{cases}\eeas where $p_1(z),\;p_2(z)$ are two non-constant polynomials.\vspace{1mm}
\par Since $D:=a_1a_4-a_2a_3\neq0$, solving the above system of equations, we get
\bea\label{e3.59} f_1(z+c)=\frac{1}{2D}\left(a_4\left(e^{p_1(z)}+e^{-p_1(z)}\right)+ia_2\left(e^{p_2(z)}-e^{-p_2(z)}\right)\right),\eea
\bea\label{e3.60} f_1(z)=\frac{1}{-2D}\left(a_3\left(e^{p_1(z)}+e^{-p_1(z)}\right)+ia_1\left(e^{p_2(z)}-e^{-p_2(z)}\right)\right),\eea
\bea\label{e3.61} f_2(z+c)=\frac{1}{2D}\left(a_4\left(e^{p_2(z)}+e^{-p_2(z)}\right)+ia_2\left(e^{p_1(z)}-e^{-p_1(z)}\right)\right)\eea and 
\bea\label{e3.62} f_2(z)=\frac{1}{-2D}\left(a_3\left(e^{p_2(z)}+e^{-p_2(z)}\right)+ia_1\left(e^{p_1(z)}-e^{-p_1(z)}\right)\right).\eea
\par From \eqref{e3.59} and \eqref{e3.60}, we obtain
\bea\label{e3.63} && \frac{ia_3}{a_1}e^{p_2(z)+p_1(z+c)}+\frac{ia_3}{a_2}e^{p_2(z)-p_1(z+c)}-\frac{a_1}{a_2}e^{p_2(z)+p_2(z+c)}+\frac{a_1}{a_2}e^{p_2(z)-p_2(z+c)}\nonumber\\&&-\frac{ia_4}{a_2}e^{p_2(z)+p_1(z)}-\frac{ia_4}{a_2}e^{p_2(z)-p_1(z)}+e^{2p_2(z)}=1.\eea
\par  From \eqref{e3.61} and \eqref{e3.62}, we obtain
\bea\label{e3.64} && \frac{ia_3}{a_1}e^{p_1(z)+p_2(z+c)}+\frac{ia_3}{a_2}e^{p_1(z)-p_2(z+c)}-\frac{a_1}{a_2}e^{p_1(z)+p_1(z+c)}+\frac{a_1}{a_2}e^{p_1(z)-p_1(z+c)}\nonumber\\&&-\frac{ia_4}{a_2}e^{p_1(z)+p_2(z)}-\frac{ia_4}{a_2}e^{p_1(z)-p_2(z)}+e^{2p_1(z)}=1.\eea
\par Now, we consider the following cases.\vspace{1mm}
\par \textbf{Case 1:} Suppose $p_2(z)-p_1(z)=k$, where $k$ is a constant in $\mathbb{C}$.\vspace{1mm}
\par Then \eqref{e3.63} and \eqref{e3.64}, respectively yield
\bea\label{e3.65}&& e^{k}\left(ia_3-a_1e^{k}\right)e^{p_1(z)+p_1(z+c)}+\left(ia_3e^{k}+a_1\right)e^{p_1(z)-p_1(z+c)}\nonumber\\&&+e^{k}\left(a_2e^{k}-ia_4\right)e^{2p_1(z)}=\left(a_2+ia_4e^{k}\right)\eea and 
\bea\label{e3.66}&& \left(ia_3e^{k}-a_1\right)e^{p_1(z)+p_1(z+c)}+\left(ia_3e^{-k}+a_1\right)e^{p_1(z)-p_1(z+c)}\nonumber\\&&+\left(a_2-ia_4e^{k}\right)e^{2p_1(z)}=\left(a_2+ia_4e^{-k}\right).\eea
\par Now, we show that all of $ia_3-a_1e^{k}$, $ia_3e^{k}+a_1$, $ a_2e^{k}-ia_4$, $a_2+ia_4e^{k}$, $ia_3e^{k}-a_1$, $ia_3e^{-k}+a_1$, $ a_2-ia_4e^{k}$ and $a_2+ia_4e^{-k}$ are non-zero.\vspace{1mm}
\par Suppose $ia_3-a_1e^{k}=0$. Then, clearly, $a_2e^{k}-ia_4=-iD/a_1\neq0$.\vspace{1mm}
\par Therefore, \eqref{e3.65} yields
\bea\label{e3.67} \left(ia_3e^{k}+a_1\right)e^{p_1(z)-p_1(z+c)}+e^{k}\left(a_2e^{k}-ia_4\right)e^{2p_1(z)}=\left(a_2+ia_4e^{k}\right).\eea
\par From \eqref{e3.67}, it is clear that $ia_3e^{k}+a_1$ is non-zero. Otherwise, $p_1(z)$ would be constant, which is not possible.\vspace{1mm}
\par Also, we claim that $a_2+ia_4e^{k}$ is non-zero. If not, then we must have from \eqref{e3.67} that \beas-\left(\frac{ia_3e^{k}+a_1}{a_2e^{k}-ia_4}\right)e^{-(p_1(z)+p_1(z+c)+k)}=1,\eeas which implies that $p_1(z)+p_1(z+c)$ and hence $p_1(z)$ is constant, which is a contradiction.\vspace{1mm}
\par Now, keeping in view of \eqref{e3.67}, we obtain \beas T\left(r,e^{p_1(z)-p_1(z+c)}\right)=T\left(r,e^{2p_1(z)}\right)+S\left(r,e^{2p_1(z)}\right).\eeas
\par Since $p_1(z)$ is a polynomial, it is easy to see that \beas N\left(r, \frac{1}{e^{p_1(z)-p_1(z+c)}}\right)=N\left(r, e^{p_1(z)-p_1(z+c)}\right)=N\left(r, \frac{1}{e^{2p_1(z)}}\right)=S\left(r,e^{p_1(z)}\right).\eeas
\par Then, keeping in view of \eqref{e3.67} and using second fundamental theorem of Nevanlinna in several complex variables, we obtain
\beas T\left(r,e^{p_1(z)-p_1(z+c)}\right)&\leq& N\left(r, \frac{1}{e^{p_1(z)-p_1(z+c)}}\right)+N\left(r, e^{p_1(z)-p_1(z+c)}\right)\\&&+N\left(r, \frac{1}{e^{p_1(z)-p_1(z+c)}-\alpha}\right)+S\left(r,e^{p_1(z)-p_1(z+c)}\right)\\&\leq& N\left(r, \frac{1}{e^{2p_1(z)}}\right)+S\left(r,e^{p_1(z)-p_1(z+c)}\right)\\&\leq& S\left(r,e^{p_1(z)-p_1(z+c)}+S\left(r,e^{2p_1(z)}\right)\right)\eeas where $\alpha=\left(a_2+ia_2e^{k}\right)/\left(ia_3e^{k}+a_1\right)$.\vspace{1mm}
\par This implies that $\left(r,e^{2p_1(z)}\right)=o\left(T\left(r,e^{2p_1(z)}\right)\right)$, which is not possible as $e^{p_1(z)}$ is transcendental entire.\vspace{1mm}
\par Hence, we conclude that $ia_3-a_1e^{k}\neq0$. Similarly, we can prove that the others are also non-zero.\vspace{1mm}
\par Now, Using Lemma \ref{lem3.1}, we obtain from \eqref{e3.65} that \bea\label{e3.68} \left(\frac{ia_3e^{k}+a_1}{a_2+ia_4e^{k}}\right)e^{p_1(z)-p_1(z+c)}=1.\eea
\par Using \eqref{e3.68} in \eqref{e3.65}, we obtain
\bea\label{e3.69} -\left(\frac{ia_3-a_1e^{k}}{a_2e^{k}-ia_4}\right)e^{-p_1(z)+p_1(z+c)}=1.\eea
\par Again, Using Lemma \ref{lem3.1}, we obtain from \eqref{e3.66} that \bea\label{e3.70} \left(\frac{ia_3e^{-k}+a_1}{a_2+ia_4e^{-k}}\right)e^{p_1(z)-p_1(z+c)}=1.\eea
\par Using \eqref{e3.70} in \eqref{e3.66}, we obtain
\bea\label{e3.71} -\left(\frac{ia_3e^{k}-a_1}{a_2-ia_4e^{k}}\right)e^{-p_1(z)+p_1(z+c)}=1.\eea
\par Keeping in view of the fact that $D\neq0$, from \eqref{e3.68} and \eqref{e3.70}, we obtain 
\beas e^{2k}=1.\eeas
\par Multiplying \eqref{e3.68} and \eqref{e3.69}, we obtain 
\beas ((a_2^2+a_4^2)-(a_1^2+a_3^2))e^{k}+i(a_2a_4-a_1a_3)(e^{2k}-1)=0.\eeas 
\par As $e^{2k}=1$, the above equation yields \beas a_2^2+a_4^2=a_1^2+a_3^2. \eeas
\par Now, in view of \eqref{e3.68}, we conclude that $p_1(z)-p_1(z+c)$ is constant. Since $p_1(z)$ is a polynomial, we may assume that \beas p_1(z)=L(z)+\Phi(t)+A,\eeas where $L(z)=\alpha_1z_1+\alpha_2z_2$, $\alpha_1,\alpha_2, A\in\mathbb{C}$ and $\Phi(t)$ is a polynomial in $t:=c_2z_1-c_1z_2$.\vspace{1mm}
\par Therefore, from \eqref{e3.68}, \eqref{e3.69}, \eqref{e3.70} and \eqref{e3.71}, we obtain 
\beas e^{L(c)}=\frac{a_1+ia_3e^{k}}{a_2+ia_4e^{k}}=\frac{a_2e^{k}-ia_4}{a_1e^{k}-ia_3}=\frac{a_1+ia_3e^{-k}}{a_2+ia_4e^{-k}}=\frac{a_2-ia_4e^{k}}{a_1-ia_3e^{k}}.\eeas
\par \textbf{Case 2:} Suppose $p_2(z)-p_1(z)$ is non-constant.\vspace{1mm}
\par We claim that $p_2(z)-p_1(z+c)$ is non-constant. If not, then let $p_2(z)-p_1(z+c)=k$, where $k$ is a constant in $\mathbb{C}$.\vspace{1mm}
\par Since $p_1(z),p_2(z)$ are non-constants polynomials, we must have $p_2(z)=L(z)+\Phi(t)+A$, where $L(z)=\alpha_1z_1+\alpha_2z_2$, $\alpha_1, \alpha_2, A\in\mathbb{C}$ and $\Phi(t)$ is a polynomial in $t:=c_2z_1-c_1z_2$. This implies that $p_1(z)=L(z)+\Phi(t)+A-L(c)-k$.\vspace{1mm}
\par Hence, $p_2(z)-p_1(z)=L(c)+k=$constant, which is a contradiction.\vspace{1mm}
\par Also, we claim that $p_2(z)-p_2(z+c)$ is constant.\vspace{1mm} 
\par Suppose on contrary, $p_2(z)-p_2(z+c)$ is non-constant. Then $p_2(z)+p_1(z+c)$ is also non-constant.\vspace{1mm}
\par Now, using Lemma \ref{lem3.1}, we obtain from \eqref{e3.63}
\bea\label{e3.72} -\frac{ia_4}{a_2}e^{p_2(z)+p_1(z)}=1.\eea 
\par Using \eqref{e3.72} in \eqref{e3.63} we obtain
\bea\label{e3.73} -ia_3e^{-k}e^{p_1(z)+p_1(z+c)}+(a_1-ia_3e^{-k})e^{p_1(z)-p_1(z+c)}&-&a_1e^{-2k}e^{-(p_1(z)+p_1(z+c))}\nonumber\\&=&a_2-ia_4e^{-k}.\eea
\par If $a_1-ia_3e^{-k}=0$, then we have from \eqref{e3.73} that \beas ia_3e^{-k}e^{2(p_1(z)+p_1(z+c))}+(a_2-ia_4e^{-k})e^{p_1(z)+p_1(z+c)}+a_1e^{-2k}=0,\eeas which implies that $p_1(z)+p_1(z+c)$ and hence $p_1(z)$ is constant, which is a contradiction. Therefore, $a_1-ia_3e^{-k}\neq0$.\vspace{1mm}

\par If $a_2-ia_4e^{-k}=0$, then \eqref{e3.73} reduces to \beas -\frac{ia_3e^{k}}{a_1}e^{2(p_1(z)+p_1(z+c))}+\frac{a_1-ia_3e^{-k}}{a_1e^{-2k}}e^{2p_1(z)}=1,\eeas which implies that $T\left(r,e^{p_1(z)+p_1(z+c)}\right)=T\left(r,e^{p_1(z)}\right)+S\left(r,e^{p_1(z)}\right)$.\vspace{1mm}
\par Now, keeping in view of the above equation, and applying second fundamental theorem of Nevanlinna in several complex variables to $e^{2(p_1(z)+p_1(z+c))}$, we obtain 
\beas T\left(r, e^{2(p_1(z)+p_1(z+c))}\right)&\leq&\ol N\left(r,e^{2(p_1(z)+p_1(z+c))}\right)+\ol N\left(r,\frac{1}{e^{2(p_1(z)+p_1(z+c))}}\right)\\&+&\ol N\left(r,\frac{1}{e^{2(p_1(z)+p_1(z+c))}-\omega}\right)+S\left(r,e^{p_1(z)+p_1(z+c)}\right)\\&\leq& \ol N\left(r,\frac{1}{e^{2p_1(z)}}\right)+S\left(r,e^{p_1(z)}\right)=S\left(r,e^{p_1(z)}\right),\eeas where $\omega=-a_1/(ia_3e^{k})$.\vspace{1mm}
\par But this implies that $T\left(r,e^{p_1(z)}\right)=S\left(r,e^{p_1(z)}\right),$ which is a contradiction. Therefore, $a_2-ia_4e^{-k}\neq0$.\vspace{1mm}

\par Now, applying Lemma \ref{lem3.1}, we obtain from \eqref{e3.73} that \beas \left(\frac{a_1-ia_3e^{-k}}{a_2-ia_4e^{-k}}\right)e^{p_1(z)-p_1(z+c)}=1,\eeas which implies that $p_1(z)-p_1(z+c)$ is constant. But, then we may assume $p_1(z)=L(z)+\Phi(t)+A$, where $L(z)$, $\Phi(t)$, $A$ are defined in Theorem \ref{t3}.\vspace{1mm}
\par Therefore, $p_2(z)=-(L(z)+\Phi(t)+A)+\text{const.}$ Hence, $p_2(z)-p_2(z+c)=L(c)$=constant, which contradicts the assumption. So, our claim is proved.\vspace{1mm}
\par Now, we consider the following subcases.\vspace{1mm}
\par \textbf{Subcase 2.1:} Suppose $p_2(z)+p_1(z)=k$, where $k$ is a constant in $\mathbb{C}$.\vspace{1mm}
\par Then from \eqref{e3.63} and \eqref{e3.64}, we obtain \bea\label{e3.74}&& \left(ia_3e^{k}+a_1\right)e^{-p_1(z)+p_1(z+c)}+e^{k}\left(ia_3-a_1e^{k}\right)e^{-(p_1(z)+p_1(z+c))}\nonumber\\&&+e^{k}\left(a_2e^{k}-ia_4\right)e^{-2p_1(z)}=\left(a_2+ia_4e^{k}\right)\eea and 
\bea\label{e3.75}&& \left(ia_3e^{k}+a_1\right)e^{p_1(z)-p_1(z+c)}+\left(ia_3e^{-k}-a_1\right)e^{p_1(z)+p_1(z+c)}\nonumber\\&&+\left(a_2-ia_4e^{-k}\right)e^{2p_1(z)}=\left(a_2+ia_4e^{-k}\right).\eea
\par In a similar manner as done in Case 1, we can prove that $ia_3e^{k}+a_1$, $ia_3-a_1e^{k}$, $a_2e^{k}-ia_4$, $a_2+ia_4e^{k}$, $ia_3e^{-k}+a_1$, $ia_3e^{-k}-a_1$, $a_2-ia_4e^{-k}$ and $a_2+ia_4e^{-k}$ are all non-constants.\vspace{1mm}
\par Now, using Lemma \ref{lem3.1}, we obtain from \eqref{e3.74} that \bea\label{e3.76} \left(\frac{ia_3e^{k}+a_1}{a_2+ia_4e^{k}}\right)e^{-p_1(z)+p_1(z+c)}=1.\eea

\par Using \eqref{e3.76} in \eqref{e3.74}, we obtain 
\bea\label{e3.77} -\left(\frac{ia_3-a_1e^{k}}{a_2e^{k}-ia_4}\right)e^{p_1(z)-p_1(z+c)}=1.\eea
\par Again, using Lemma \ref{lem3.1}, we obtain from \eqref{e3.75} that \bea\label{e3.78} \left(\frac{ia_3e^{-k}+a_1}{a_2+ia_4e^{k}}\right)e^{p_1(z)-p_1(z+c)}=1.\eea

\par Using \eqref{e3.78} in \eqref{e3.75}, we obtain 
\bea\label{e3.79} -\left(\frac{ia_3e^{-k}-a_1}{a_2-ia_4e^{-k}}\right)e^{-p_1(z)+p_1(z+c)}=1.\eea
\par Now, from \eqref{e3.76} and \eqref{e3.79}, we obtain $e^{2k}=1$. \vspace{1mm}
\par From \eqref{e3.76}--\eqref{e3.79}, we obtain
\bea\label{e3.80} \frac{a_2+ia_4e^{k}}{a_1+ia_3e^{k}}=\frac{a_1e^{k}-ia_3}{a_2e^{k}-ia_4}=\frac{ia_3e^{k}+a_1}{a_2+ia_4e^{k}}=\frac{a_2-ia_4e^{-k}}{a_1-ia_3e^{-k}}.\eea
\par Now, if $e^{k}=1$, then from second and fifth terms of \eqref{e3.80}, we get $D=0$, which is a contradiction.\vspace{1mm}
\par Similarly, for $e^{k}=-1$, we get a contradiction.\vspace{1mm}
\par \textbf{Subcase 2.2:} Suppose $p_2(z)+p_1(z)$ is non-constant. As we already proved that $p_2(z)-p_2(z+c)$ is constant, let it be $k$, $k\in\mathbb{C}$.\vspace{1mm}
\par Then \eqref{e3.63} yields
\beas ia_3e^{p_2(z)+p_1(z+c)}&+&ia_3e^{p_2(z)-p_1(z+c)}-a_1e^{p_2(z)+p_2(z+c)}-ia_4e^{p_2(z)+p_1(z)}\\&-&ia_4e^{p_2(z)-p_1(z)}+a_2e^{2p_2(z)}=a_2+a_1e^{k}.\eeas

\par Using Lemma \ref{lem3.1}, we get from the above equation that \beas \left(\frac{ia_3}{a_2+a_1e^{k}}\right)e^{p_2(z)+p_1(z+c)}=1,\eeas which implies that $p_2(z)+p_1(z+c)$ is constant. But the $p_2(z)+p_1(z)$ must be constant, which is a contradiction to our assumption.\vspace{1mm}
\par \textbf{Subcase 2.3:} Suppose $p_2(z)+p_1(z+c)=k$, where $k$ is a constant in $\mathbb{C}$. Then we may write $p_2(z)=L(z)+\Phi(t)+A$, where $L(z)$, $\Phi(t)$, $A$ are defined in Theorem \ref{t3}.\vspace{1mm}
\par Therefore, it may be easily seen that $p_1(z)=-(L(z)+\Phi(t)+A)+L(c)+k$, $p_2(z)+p_1(z)=L(c)+k$ and $p_2(z)-p_2(z+c)=-L(c)$. Thus, we can obtain from \eqref{e3.63} that \beas \alpha e^{2p_2(z)}=\beta,\eeas which implies that $p_2(z)$ is constant, a contradiction.\vspace{1mm}
\par \textbf{Subcase 2.4:} Suppose $p_2(z)+p_1(z+c)$ is non-constant. 
As we know that $p_2(z)-p_2(z+c)$ is constant, let it be $k\in\mathbb{C}$.\vspace{1mm}
\par Therefore, \eqref{e3.63} reduces to
\beas ia_3e^{p_2(z)+p_1(z+c)}+ia_3e^{p_2(z)-p_1(z+c)}-a_1e^{p_2(z)+p_2(z+c)}-ia_4e^{p_2(z)+p_1(z)}\\-ia_4e^{p_2(z)-p_1(z)}+a_2e^{2p_2(z)}=a_2-a_1e^{k},\eeas where one can easily verify that $a_2-a_1e^{k}$ is non-zero.\vspace{1mm}
\par Using Lemma \ref{lem3.1}, we obtain from the above equation 
\beas \left(\frac{-ia_4}{a_2-a_1e^k}\right)e^{p_2(z)+p_1(z)}=1,\eeas which implies that $p_2(z)+p_1(z)=$constant.\vspace{1mm}
\par since $p_2(z)-p_2(z+c)=k$, we may assume that $p_2(z)=L(z)+\Phi(t)+A$, where $L(z)$, $\Phi(t)$ and $A$ are defined in Theorem \ref{t3}. But then we can get that $p_2(z)+p_1(z+c)$ is constant, which contradicts our assumption.\vspace{1mm}
\par This completes the proof of the theorem.\vspace{1mm}
\begin{proof}[\textbf{Proof of Theorem $\ref{t4}$}]
The proof of this theorem can be carried out with similar arguments as in the proof of Theorem 1.1 of \cite{Xu-Liu-Li-JMAA-2020}. So, we omit the details. 
\end{proof}
\end{proof}
\vspace{1.5mm}

\noindent{\bf Acknowledgment:} The authors would like to thank the referee(s) for the helpful suggestions and comments to improve the exposition of the paper.

\end{document}